\documentclass[12pt,a4paper]{article} 
\usepackage{amsfonts,amssymb,amsmath,amscd,stmaryrd,latexsym,makeidx,theorem}
\usepackage{amsmath,pstricks,latexsym,theorem,epsfig} 
\usepackage{amsfonts,amssymb,latexsym,epsfig,amsmath}
\usepackage{amsmath,latexsym,theorem} 
\usepackage{graphicx}

\usepackage{amsmath,latexsym,theorem} 
\newtheorem{Theorem}{Theorem}[section] 
\newtheorem{Definition}{Definition}[section] 

\newtheorem{Proposition}{Proposition}[section] 
\newtheorem{Lemma}{Lemma}[section] 
\newtheorem{Corollary}{Corollary}[section] 
 
\theorembodyfont{\rmfamily} 
\newtheorem{Remark}{Remark}[section]

\newcommand{\rec}[1]{{(\ref{#1})}} 
\newcommand{\lapfr}{(-\Delta)^\frac{1}{2}} 
\newcommand{\ba}{\begin{array}} 
\newcommand{\ea}{\end{array}}

 \def\11{1\!\!1}

\newcommand{\be}{\begin{equation}}
\newcommand{\ee}{\end{equation}}
\newcommand{\bes}{\begin{equation*}}
\newcommand{\ees}{\end{equation*}}
 \newcommand{\R}{\mathbb{R}}

\newcommand{\C}{\mathbb{C}}
\newcommand{\Z}{\mathbb{Z}}

\def\ti{\tilde}
\def\lf{\left}
\def\rg{\right}

\def\al{\alpha}

\def\ds{\displaystyle}

\def\Om{\Omega}

\def\p{\partial}

\def\ba{\vec{a}}

\begin{document} 
 
 \title{\bf  A Pohozaev-type formula and Quantization of Horizontal Half-Harmonic Maps}
\author{ Francesca Da Lio\thanks{Department of Mathematics, ETH Z\"urich, R\"amistrasse 101, 8092 Z\"urich, Switzerland.}  \and Paul Laurain \thanks{Institut de Math\'ematiques de Jussieu, Paris 7,75205 PARIS Cedex 13
France} \and Tristan Riviere$^*$   }
\maketitle 
 
\begin{abstract}
 In a recent paper \cite{DLR3}  the first and the third  authors introduced the notion of {\em horizontal $\alpha$ harmonic map}, $\alpha\ge 1/2$ with respect to a given $C^1$ planes distribution $P_T$ on all ${\mathbb R}^m$.  These are maps $u\in \dot H^{\alpha}({\mathbb R}^k,{\R}^m)$, $\alpha\ge 1/2$, satisfying 
 $
P_T\nabla u=\nabla u$ and $P_T(u)(-\Delta)^{\alpha}u=0$ in $ {\mathcal D}'({{\mathbb R}}^k).$
The goal of this paper is to investigate compactness and quantization properties of sequences of horizontal $1/2$ harmonic maps $u_k $ in $1$D such that 
 $\|u_k\|_{\dot H^{1/2}(\R)}\le C$ and $\|(-\Delta)^{1/2}u_k\|_{L^1(\R)}\le C\,.$ We show that there exist a horizontal $1/2$ harmonic map $u_\infty$
 and a possibly empty set $\{a_1,\ldots, a_\ell\}$, $\ell\ge 1\,,$ such that up to subsequence
\begin{equation} 
u_{k }\to u_\infty \quad \mbox{in $\dot W^{1/2,p}_{loc}(\R\setminus\{a_1,\ldots, a_{\ell}\}),~~p\ge 2$ as $k \to +\infty.$}
\end{equation}
 Moreover 
there is a family $ \tilde u_{\infty}^{i,j} $ of horizontal $1/2$-harmonic maps $ ( i\in\{1,\ldots,\ell\}, j\in\{1,\ldots ,N_i\}), $ such that up to subsequence
$$
\left\| (-\Delta)^{1/4}\left(u_k-u_{\infty}-\sum_{i,j} \tilde u_{\infty}^{i,j}((x-x^k_{i,j})/r^k_{i,j})\right)\right\|_{L^2_{loc}(\R)}\to 0,~~\mbox{as $k \to +\infty$}\,
$$
for some sequences $r^k_{i,j}\rightarrow 0$ and $x^k_{i,j}\to a_i$ as $k\to \infty.$\par
The quantization analysis is obtained through a precise asymptotic development  of the energy of $u_k$  in the neck region and a subtle application of new  Pohozaev-type formulae.
 \end{abstract}

{\noindent
 }\medskip
 \tableofcontents
 \section{Introduction}
 In a recent paper \cite{DLR3}  the first and the third  authors introduced the notion of {\em horizontal $\alpha$ harmonic map}, $\alpha\ge 1/2$ with respect to a given $C^1$ planes distribution. 
Precisely  we consider  $P_T\in C^1({\R}^m,M_m({\R}))$ and $P_N\in C^1({\R}^m,M_m({\R}))$ such that
\be
\label{I.1}\lf\{
\begin{array}{l}
P_T\circ P_T=P_T\quad P_N\circ P_N=P_N\\[3mm]
P_T+P_N=I_m\\[3mm]
\forall\, z\in {\R}^m\quad\forall\, U,V\in T_z{\R}^m\quad <P_TU,P_NV>=0\\[3mm]
\|\p_z P_T\|_{L^\infty({\R}^m)}<+\infty
\end{array}
\rg.
\ee
where $<\cdot,\cdot>$ denotes the standard scalar product in ${\R}^m$. In other words $P_T$ is a $C^1$ map into the orthogonal projections of ${\R}^m$.  For such a distribution of projections $P_T$
we denote by
\[
n:=\mbox{rank}(P_T).
\]
Such a distribution identifies naturally with the distribution of $n-$planes given by the images of $P_T$ (or the Kernel of $P_T$) and conversely,
any $C^1$ distribution of $n-$dimensional planes defines uniquely $P_T$ satisfying (\ref{I.1}).

For any $\al\ge 1/2$ and for $k\ge 1$ we define the space of {\bf $H^\al$-Sobolev horizontal maps }
\[
{\frak{H}}^{\al} ({\R}^k):=\lf\{u\in H^{\al}({\R}^k,{\R}^m)\quad;\quad P_N(u)\nabla u=0\quad\mbox{ in }{\mathcal D}'({\R}^k)\rg\}
\]
 Observe that this definition makes sense since we have respectively $P_N\circ u\in H^{\al}(\R^k,M_m({\R}))$ and ${\nabla u}\in H^{\al-1}({\R}^k,{\R}^m)$.
 \begin{Definition}
 \label{df-legend-harm}
 Given a $C^1$ plane distribution $P_T$ in ${\R}^m$ satisfying (\ref{I.1}), a map $u$ in the space ${\frak H}^{\al}({\R}^k)$ is called  {\bf {horizontal  $\al$-harmonic} }with respect to $P_T$
 if
 \be
 \label{int-I.6}
 \forall\, i=1\cdots m\quad\quad\sum_{j=1}^mP_T^{ij}(u)(-\Delta)^\al u_j=0\quad\quad\mbox{in }{\mathcal D}'({\R}^k)
 \ee
and we shall use the following notation
\[
 P_T(u)\,(-\Delta)^\al u=0\quad\quad\mbox{in }{\mathcal D}'({\R}^k).\hfill \Box
\] 
 \end{Definition}
 When the plane distribution $P_T$ is {\it integrable} that is to say when
\be
\label{I.2}
\forall \ X,Y\in C^1({\R}^m,{\R}^m)\quad P_N [P_T\, X,P_T\,Y]\equiv 0
\ee
where $[\cdot,\cdot]$ denotes the Lie Bracket of vector-fields, by using Fr\"obenius theorem the planes distribution corresponds to
the tangent plane  distribution of a $n-$dimensional {\it foliation } ${\mathcal F}$, (see e.g \cite{Lang}). A smooth map $u$ in ${\frak H}^{\al}({\R}^m)$ takes values  everywhere
into a {\it leaf } of ${\mathcal F}$ that we denote $N^n$ and we are back to the classical theory of $\al$ harmonic maps  into
manifolds. We recall that the notion of weak $1/2$ harmonic maps into a $n$-dimensional closed manifolds ${N}^n\subset \R^m$ has been introduced by the first and third author in
\cite{DLR1,DLR2}.  
These maps
are critical points of the fractional energy on ${\R}^k$
\begin{equation}\label{fracenergy}
E^{1/2}(u):=\int_{{\R}^k}|(-\Delta)^{1/4}u|^2\ dx^k
\end{equation}
within
\[
H^{1/2}({\R}^k,N^n):=\lf\{ u\in H^{1/2}({\R}^k,{\R}^m)\ ;\ u(x)\in N^n\ \mbox{ for a. e. }x\in {\R}^k\rg\}.
\]
The corresponding Euler-Lagrange equation is given by
\be
\label{int-I.3}
\nu(u)\wedge(-\Delta)^{1/2} u=0\quad \quad\mbox{in }{\mathcal D}'({\R}^k),
\ee
 where $\nu(z)$ is the Gauss Maps at $z\in\cal{N}$ taking values into the grassmannian $\ti{G}r_{m-n}({\R}^m)$ of oriented $m-n$ planes in ${\R}^m$ which
 is given by the oriented normal $m-n-$plane to $T_z\cal{N}\,.$
One of the main results obtained in \cite{DLR3} is the following
\begin{Theorem}
\label{th-I.1}
Let $P_T$ be a $C^1$ distribution of planes (or projections) satisfying (\ref{I.1}). Any map $u\in {\frak H}^{1/2}({\R})$ (resp.  $u\in {\frak H}^{1}({\R^2})$)  satisfying
\be
\label{I.4}
P_T(u)\,(-\Delta)^{1/2}u=0\quad\mbox{ in }{\mathcal D}'({\R})
\ee 
(resp.
\be
\label{I.5}
P_T(u)\,\Delta u=0  \quad\mbox{ in }{\mathcal D}'({\R^2}))
\ee 
is in $\cap_{\delta<1}C^{0,\delta}({\R})$, (resp. $\cap_{\delta<1}C^{0,\delta}({\R^2})$). \hfill $\Box$
\end{Theorem}
In order to prove Theorem \ref{th-I.1} the authors in \cite{DLR3} use the following two key properties satisfied by respectively horizontal harmonic and horizontal $1/2$-harmonic maps.\par
Horizontal harmonic maps   satisfy an elliptic Schr\"odinger type system with an antisymmetric potential $\Om\in L^{2}({\R}^k,{\R}^k\otimes so(m))$ of the form
\be
\label{int-I.7}
-\Delta u=\Omega(P_T)\cdot\nabla u.
\ee
Hence, following the analysis in \cite{Riv} the authors  deduced in two dimension  the local existence on a disc $D^2$ of $A(P_T)\in L^\infty\cap W^{1,2}(D^2,Gl_m({\R}))$ and $B(P_T)\in W^{1,2}(D^2,M_m({\R}))$ such that
\be
\label{int-I.8}
\mbox{div}\lf(A(P_T)\,\nabla u\rg)=\nabla^\perp B(P_T)\cdot\nabla u
\ee
from which the regularity of $u$ can be deduced using Wente's {\it Integrability by compensation} which can be summarized in the following estimate
\be
\label{int-I.9}
\|\nabla^\perp B\cdot\nabla  u\|_{H^{-1}(D^2)}\le \, C\, \|\nabla B\|_{L^2(D^2)}\ \|\nabla u\|_{L^2(D^2)}.
\ee
A similar property is satisfied by horizontal $1/2$-harmonic maps. Precisely in \cite{DLR3}   conservation laws corresponding to (\ref{int-I.8}) but for general {\it horizontal $1/2-$harmonic maps} have been discovered: locally, modulo some smoother terms 
coming from the application of non-local operators on cut-off functions, the authors construct $A(P_T)\in L^{\infty}\cap H^{1/2}(\R,Gl_m({\R}))$ and $B(P_T)\in H^{1/2}(\R,M_m({\R}))$ such that
\be
\label{int-I.9}
(-\Delta)^{1/4}(A(P_T)\, v)=\mathcal{J}(B(P_T),v) +\mbox{cut-off},
\ee
where $v:=(P_T\,(-\Delta)^{1/4}v,{\mathcal R}\, (P_N(-\Delta)^{1/4}v))$ and ${\mathcal R}$ denotes the Riesz operator and ${\mathcal J}$ is a bilinear pseudo-differential operator satisfying
\be
\label{int-I.10}
\|\mathcal{J}(B,v)\|_{H^{-1/2}({\R})}\le C\, \|(-\Delta)^{1/4} B\|_{L^2({\R})}\, \|v\|_{L^2({\R})}.
\ee

Moreover  by assuming that $P_T\in C^2(\R^m)$and $\|\p_zP_T\|_{ L^{\infty}(\R)}<+\infty$ and by bootstrapping  arguments  one gets that every horizontal $1/2$ harmonic map $u\in {\frak H}^{1/2}({\R})$ is $C^{1,\alpha}_{loc}(\R)$, for every $\alpha<1$, (see \cite{DL3}).\par
We also remark that if   $\Pi\colon S^1\setminus \{-i\}\to \R$, $\Pi(\cos(\theta)+i \sin(\theta))=\frac{\cos(\theta)}{1+\sin(\theta)}$ is the classical stereographic projection whose inverse is given by
 \begin{equation}\label{formulaPi}
\Pi^{-1}(x)=\frac{2x}{1+x^2}+i\left(-1+\frac{2}{1+x^2}\right).
\end{equation}
then the  following relation between the $1/2$ Laplacian in $\R$ and in $S^1$ holds:
 \begin{Proposition}[Proposition 4.1, \cite{DLMR}] \label{proppull} Given $u:\R\to\R^m$ set $v:=u\circ \Pi:S^1\to\R^m$. Then $u\in L_{\frac{1}{2}}(\R)$\footnote{We recall that 
 $L_\frac{1}{2}(\R):=\left\{u\in L^1_{loc}(\R):\int_{\R}\frac{|u(x)|}{1+x^2}dx<\infty   \right\}$} if and only if $v\in L^1(S^1)$. In this case
\begin{equation}\label{eqlapv1}
\lapfr_{S^1} v(e^{i\theta})=\frac{((-\Delta)_{\R}^\frac12u)(\Pi(e^{i\theta}))}{1+\sin\theta}  ~~\mbox{in $\mathcal{D}'(S^1\setminus\{-i\})$},
\end{equation}
Observe that $1+\sin(\theta)=|\Pi^{\prime}(\theta)|,$ and hence we have
$$\int_{S^1} \lapfr v(e^{i\theta})\,\varphi(e^{i\theta})\ d\theta=\int_{\R}\lapfr u(x)\ \varphi\circ\Pi^{-1}(x)\ dx\quad \text{for every }\varphi\in C^\infty_0(S^1\setminus\{-i\}).$$\end{Proposition}
From Proposition \ref{proppull}  it follows that $u\in {\frak H}^{1/2}({\R})$ is a horizontal $1/2$ harmonic map in $\R$  if and only if  $v:=u\circ \Pi\in {\frak H}^{1/2}({S^1})$  is a horizontal $1/2$ harmonic map in $S^1$.
\par
\medskip
The goal of this paper is to investigate compactness and quantization properties of sequences of horizontal $1/2$ harmonic maps $u_k\in {\frak{H}}^{1/2} ({\R}).$ Our  main result is the following:
\begin{Theorem}
\label{th-I.2}
Let $u_k\in {\frak{H}}^{1/2} ({\R})$ be a sequence of horizontal $1/2$-harmonic maps such that 
\be\label{hypcom}\|u_k\|_{\dot H^{1/2}}\le C,~~~\|(-\Delta)^{1/2}u_k\|_{L^1}\le C\,.\ee
 Then it holds:\par
\begin{enumerate}
\item There exist $u_{\infty}\in {\frak{H}}^{1/2} ({\R})$ and a possibly empty set $\{a_1,\ldots, a_\ell\}$, $\ell\ge 1\,,$ such that up to subsequence
\begin{equation}\label{conv}
u_{k }\to u_\infty \quad \mbox{in $\dot W^{1/2,p}_{loc}(\R\setminus\{a_1,\ldots, a_{\ell}\}),~~p\ge 2$ as $k \to +\infty$}
\end{equation}
and 
\begin{equation}\label{provv1}
P_T(u_{\infty})(-\Delta)^{1/2}  u_{\infty}=0,\quad \mbox{in ${\cal{D}}^{\prime}(\R )$\,.}
\end{equation}
\item There is a family $ \tilde u_{\infty}^{i,j}\in \dot {\frak{H}}^{1/2} ({\R})$ of horizontal $1/2$-harmonic maps $ ( i\in\{1,\ldots,\ell\}, j\in\{1,\ldots ,N_i\}), $ such that up to subsequence
\begin{equation}\label{finalquantintr}
\left\| (-\Delta)^{1/4}\left(u_k-u_{\infty}-\sum_{i,j} \tilde u_{\infty}^{i,j}((x-x^k_{i,j})/r^k_{i,j})\right)\right\|_{L^2_{loc}(\R)}\to 0,~~\mbox{as $k \to +\infty$}\,.
\end{equation}
for some sequences $r^k_{i,j}\rightarrow 0$ and $x^k_{i,j}\in {\R}$.
\end{enumerate}
\end{Theorem}
We would like to make some comments and remarks on Theorem \ref{th-I.2}.\par
 {\bf 1.} We first mention that the condition $\|(-\Delta)^{1/2}u_k\|_{L^1}\le C$ is always satisfied in the case the maps $u_k$ take values into a closed manifold of $\R^m$ (case of sequences of $1/2$ harmonic maps) as soon as $\|u_k\|_{\dot H^{1/2}}\le C.$ This follows from the fact that if $u$ is a $1/2$-harmonic maps with values into a closed manifold of ${\cal{N}}^n$ of $\R^m$ then
the following inequality holds (see {Proposition} \ref{prL1})
\begin{equation}
\label{L^1}
\|(-\Delta)^{1/2}u\|_{L^1(\R)}\le C\|(-\Delta)^{1/4}u\|^2_{L^2(\R)}\quad.
\end{equation}
Hence we have the following corollary
\begin{Corollary}
Let ${\mathcal N}^n$ be a closed $C^2$ submanifold of ${\R}^m$ and let Let $u_k\in {{H}}^{1/2} ({\R},{\mathcal N}^n)$ be a sequence of  l $1/2$-harmonic maps such that 
\be\label{hypcom-1}\|u_k\|_{\dot H^{1/2}}\le C\ee
 then the conclusions of theorem~\ref{th-I.2} hold. In particular modulo extraction of a subsequence  we have the following energy identity:  \begin{equation}
\lim_{k\rightarrow +\infty}\int_{\R} |(-\Delta)^{1/4}u_k|^2\ dx=\int_{\R} |(-\Delta)^{1/4}u_\infty|^2\ dx+\sum_{i,j}\int_{\R} |(-\Delta)^{1/4}\tilde u_\infty^{i,j}|^2\ dx
\end{equation}
where $\tilde u_\infty^{i,j}$ are the {\bf bubbles} associated to the weak convergence.
\end{Corollary}

For the moment it remains open to know whether the bound (\ref{L^1}) holds or not in the general case of  horizontal $1/2$-harmonic maps.\par
{\bf 2.} The compactness issue (first part of Theorem \ref{th-I.2}) is quite standard. The most delicate part is the quantization analysis consisting in verifying that there is no dissipation of the energy in the region between $u_{\infty}$ and
the {\em bubbles} $\tilde u_{\infty}^{i,j}$ and between the bubbles themselves (the so-called {\em neck-regions}). The strategy of the proof of theorem~\ref{th-I.2} is presented in the next section. One important tool we are using for proving theorem~\ref{th-I.2} is a new {\bf Pohozaev identity} for the half Laplacian in 1 dimension.
\begin{Theorem}\label{th-I.4}{\bf [Pohozaev Identity in $\R$]}
Let $u\in W^{1,2}(\R,\R^m)$ be  such that
\begin{equation}\label{Iharmequation1}
\frac{du}{dx}\cdot (-\Delta)^{1/2} u=0~~\mbox{a.e in $\R$.}
\end{equation}
Assume that
\begin{equation}\label{IcondPoh}
\int_{\R}|u-u_0 |dx <+\infty,~~\int_{\R}\left\vert \frac{du}{dx}(x)\right\vert\, dx <+\infty
\end{equation}
Then the following identity holds
\begin{equation}\label{I-identityPR}
 \left|\int_{x\in\R} \frac{x^2-t^2}{(x^2+t^2)^2}u(x)dx\right|^2=\left|\int_{x\in\R}  \frac{2xt}{(x^2+t^2)^2}  u(x)dx\right|^2.\end{equation}
\end{Theorem}
  We observe that the conditions \rec{IcondPoh} are satisfied by the $1/2$-harmonic maps with valued into a closed $C^2 $ sub-manifold.
By means of the stereographic projection we get an analogous formula in $S^1$.
\begin{Theorem}
\label{Phohozaev-S^1}
{\bf [Pohozaev Identity on $S^1$]}
Let $u$ be a $W^{1,2}$ map from $S^1$ into ${\R}^m$ satisfying
\be
\label{pert}
\frac{d u}{d\theta}\cdot (-\Delta)^{1/2} u=0\quad \mbox{ a. e. on }S^1
\ee
then the following identity holds
\be
\label{poho-s1}
\left|\int_0^{2\pi}u(\theta)\, \cos\theta\ d\theta\right|^2=\left|\int_0^{2\pi}u(\theta)\, \sin\theta\ d\theta\right|^2
\ee
\hfill $\Box$
\end{Theorem}
We have now to give some explanations why these identities belong to the {\it Pohozaev identities} family.  These identities are produced by the conformal invariance of the highest
order derivative term in the Lagrangian from which the Euler Lagrange is issued : in $2$D for instance
\[
E(u)=\int_{\R^2} |\nabla u|^2 dx^2
\] 
is conformal invariant, whereas 
\[
E^{1/2}(u)=\int_{\R} |(-\Delta)^{1/4}u|^2\ dx
\]
is conformal invariant  in $1$D. The infinitesimal perturbations issued from the dilations produce the following infinitesimal variations of these highest order terms respectively
\[
\sum_{i=1}^2 x_i\,\frac{\p u}{\p x_i}\cdot\Delta u\quad\mbox{ in $2$D}\quad\mbox{ and }\quad x\frac{d u}{dx}\cdot(-\Delta)^{1/2}u\quad\mbox{ in  $1$D}
\]
Being a critical point respectively of $E$ in 2 D and $E^{1/2}$ in $1D$ and assuming enough regularity gives respectively
\[
\sum_{i=1}^2 x_i\,\frac{\p u}{\p x_i}\cdot\Delta u=0\quad\mbox{ in $2$D}\quad\mbox{ and }\quad x\frac{d u}{dx}\cdot(-\Delta)^{1/2}u=0\quad\mbox{ in  $1$D}
\]
In two dimensions, integrating this identity on a ball $B(x_0,r)$ gives the following {\bf balancing law} between the {\bf radial} part and the {\bf angular} part of the energy
classically known as {\bf Pohozaev identity}.
\begin{Theorem}\label{th-I.6} 
Let $u\in W^{2,2} ({\R}^2,\R^m)$ such that
\begin{equation}\label{IharmequationR2}
\sum_{i=1}^2 x_i\,\frac{\partial u}{\partial x_i }\cdot \Delta u=0~~\mbox{a.e in $B(0,1)$.}
\end{equation}
 Then it holds
\begin{equation}
\label{PI}
\int_{\partial B(x_0,r)} \left\vert \frac{1}{r}\frac{\partial u}{\partial\theta}\right\vert^2 d\theta=\int_{\partial B(x_0,r)}\left\vert\frac{\partial u}{\partial r}\right\vert^2 d\theta
\ee
for all $r\in [0,1].$
\end{Theorem}
In 1 dimension one might wonder what corresponds to the $2$ dimensional dichotomy between {\bf radial} and {\bf angular} parts. We  illustrate below  the  correspondence
of dichotomies respectively in $1$ and $2$ dimensions.

\[
\begin{array}{ccc}
\mbox{$2$D}& \longleftrightarrow & \mbox{ $1$D}\\[3mm]
\ds\mbox{ {\bf radial} : } \quad\frac{\p u}{\p r}& \longleftrightarrow & \mbox{ {\bf symmetric} part of }u\quad:\quad u^+(x):= \frac{u(x)+u(-x)}{2}\\[3mm]
\ds\mbox{{\bf  angular} : }\quad \frac{\p u}{\p \theta}& \longleftrightarrow & \mbox{ {\bf antisymmetric} part of }u\quad:\quad  u^-(x):= \frac{u(x)-u(-x)}{2}
\end{array}
\]
 Observe moreover that our {\bf Pohozaev identity in $1$D} (\ref{poho-s1}) can be rewritten as a {\bf balancing law} between the {\bf symmetric} part and the {\bf antisymmetric} part
of $u$.
\be
\label{poho-s2}
\left|\int_0^{2\pi}u^+(\theta)\, \cos\theta\ d\theta\right|^2=\left|\int_0^{2\pi}u^-(\theta)\, \sin\theta\ d\theta\right|^2
\ee
This law {\underbar is not invariant} under the action of the M\"obius group but the condition (\ref{pert}) is. Applying for instance rotations by an arbitrary angle $\al\in{\R}$, the identity (\ref{poho-s1})
implies
\be
\label{poho-s2}
\left\{
\begin{array}{l}
\ds |u_1|=|u_{-1}|\\[3mm]
\ds u_1\cdot u_{-1}=0
\end{array}
\right.
\ee
where 
\[
\left\{\begin{array}{l}
\ds u_1:=\frac{1}{2\pi}\int_0^{2\pi}u(\theta)\, \cos\theta\ d\theta\\[3mm]
\ds u_{-1}=\frac{1}{2\pi}\int_0^{2\pi}u(\theta)\, \sin\theta\ d\theta
\end{array}
\right.
\]
The previous implies that there are ``as many'' Pohozaev identities
as elements in this group minus the action of rotations, that is there are a 3-1=2=D+1 dimensional family of identities exactly as in the 2-D case where there are
exactly as many {\bf Pohozaev identities} (\ref{PI}) as choices of center $x_0\in {\R}^2$ and radius $r>0$ (which is again a $D+1=3$ dimensional space).

\subsection{The strategy of the proof of theorem~\ref{th-I.2}.}
\label{strategy}

We first recall the definitions of a bubble and a neck region.
 
 \begin{Definition}[Bubble]
 A {\bf Bubble} is a {\bf non-constant}   horizontal  $1/2$-harmonic map $u\in {\frak{H}}^{1/2} ({\R})$\,.   \end{Definition}
 \begin{Definition}[Neck region]
 A neck region for a  sequence  $f_k\in L^2(\R)$ is the union of finite   degenerate annuli of the type $A_k(x)=B(x,R_k)\setminus B(x,r_k)$ with $r_k\to 0$ and $\frac{R_k}{r_k}\to +\infty$ as $k \to +\infty$ such that  
 \begin{equation}
 \label{neck}
\lim_{\Lambda\to \infty}\lim_{k\to \infty}\left( \sup_{\rho\in[\Lambda r_k,(2\Lambda)^{-1}R_k]}\int_{B(x,2\rho)\setminus B(x,\rho)}|f_k|^2 dx\right)^{1/2}=0\,.
 \end{equation}
 \end{Definition}

 \medskip
The main achievement of the present work is to show, under the assumptions of theorem~\ref{th-I.2}, that (\ref{neck}) for $f_k:=(-\Delta)^{1/4}u_k$ can be improved to
\begin{equation}\label{limitneck}
\lim_{\Lambda\to \infty}\lim_{k\to \infty}\left(\int_{B(x,\frac{R_k}{\Lambda})\setminus B(x,\Lambda r_r)}|f_k|^2 dx\right)^{1/2}=0.
\end{equation}\par
The proof of estimate of  \rec{limitneck} will be  the aim of Section \ref{cq}.\par
{\bf 3.}  Theorem \ref{th-I.2} has been proved by the first author in \cite{DL2} in the case of sequences of  $1/2$-harmonic maps $u_k$  with values into  the ${\cal{S}}^{m-1}$ sphere, see also \cite{LP}. In this case in order to prove
\rec{limitneck} we use the duality  of the Lorentz spaces $L^{2,1}-L^{2,\infty}.$\footnote{ see section \ref{secdef} for a definition.}

We first show that the $L^{2,\infty}$ norm of the $u_k$ is arbitrary small in the  {\em neck region} 
and then we use the fact that  $ L^{2,1}$ norm of $1/2$-harmonic maps    with values into a sphere is uniformly globally bounded. This last global estimate follows directly by
the formulation of the $1/2$-harmonic maps equation that the first and third author discovered in \cite{DLR1} in terms of special algebraic quantities (three-commutators) satisfying  particular integrability compensation properties. \par
 
The fact that the 
$L^{2,\infty}$ norm of a sequence of horizontal $1/2$-harmonic maps $u_k$ is arbitrary small in {\em neck regions}  still holds in the general case, precisely we have:
 \begin{Lemma}[$L^{2,\infty}$ estimates]\label{neck2}
   There exists  $ \delta >0$ such that for any sequence of horizontal  $1/2$-harmonic maps $u_k\in{\frak{H}}^{1/2} ({\R})$ satisfying 
     $$ \sup_{\rho\in[\Lambda r_k,{(2\Lambda)}^{-1}R_k]}\left( \int_{B(0,2\rho)\setminus B(0,\rho)}| 
 (-\Delta)^{1/4} u_k|^2dx\right)^{1/2}<\delta $$
with  $r_k\to 0$ and $\frac{R_k}{r_k}\to +\infty$ as $k \to +\infty$ and for all $\Lambda>1$ with  $\Lambda r_k<(2\Lambda)^{-1}R_k$
then
\begin{equation}\label{linftyest}
\limsup_{\Lambda\to\infty}\limsup_{k\to\infty } \| (-\Delta)^{1/4} u\|_{L^{2,\infty}(B(0,(2\Lambda)^{-1}R_k)\setminus B(0,\Lambda r_k)} =0.\end{equation}
   \end{Lemma}
Unfortunately  we do not know if  such a   a global  uniform $ L^{2,1}$ estimate exists even in the case
of $1/2$-harmonic maps with values into a closed submanifold of $\R^m.$ 
 
To overcome  of this  lack of a global $L^{2,1}$ estimate, we show a precise asymptotic development  of $(-\Delta)^{1/4}u_k$ in annuli
$A_{r_k,R_k}=B(x,R_k)\setminus B(x,r_k)$ where the $L^2$ norm of $(-\Delta)^{1/4}u_k$ is small.
To get such an asymptotic development  we argue as follows. We first prove the following result (for simplicity we will consider an annulus centered at $x=0$).
\begin{Theorem}\label{th-I.3}
Let $u_k\in{\frak{H}}^{1/2} ({\R})$ be a sequence of horizontal $1/2$-harmonic maps satisfying the assumptions of Theorem \ref{th-I.1}. There exists $\delta>0$ such that if     
 \begin{eqnarray}\label{quant0}
      \|(-\Delta)^{1/4}u_k\|_{L^2(B(0,R_k)\setminus B(0,r_k)}&<&\delta,
     \end{eqnarray}
   with $r_k\to 0$ and $\frac{R_k}{r_k}\to +\infty$ as $k \to +\infty$ , then for all $\Lambda>1$ with  $\Lambda r_k<(2\Lambda)^{-1}R_k$  and $x\in B(0,(2\Lambda)^{-1}R_k)\setminus B(0,\Lambda r_k)$ we have
      \begin{equation}\label{devel}
   (-\Delta)^{1/4}u_k(x)= \frac{a^+_k(x)\overrightarrow{c_{r_k}}}{|x|^{1/2}} +  h_k+g_k,
    \end{equation}
where $g_k\in L^2(\R)$ with
$\displaystyle \limsup_{k\to\infty}\|g_k\|_{L^2} =0,$ 
\be\label{esckI}
\displaystyle \overrightarrow{c_{r_k}}=O\left(\left(\log\left(\frac{R_k}{2\Lambda^2r_k}\right)\right)^{-1/2}\right),~~\mbox{as}~~ k\to +\infty,\Lambda\to +\infty\ee $a^+_k\in L^{\infty}\cap \dot H^{1/2}(\R,Gl_m({\R}))$, $h_k\in L^{2,1}(A_{\Lambda r_k,(2\Lambda)^{-1}R_k})$,
\be\limsup_{\Lambda\to\infty}\limsup_{k\to+\infty} \|h_k\|_{L^{2,1}(A_{\Lambda r_k,(2\Lambda)^{-1}R_k})}<\infty\ee
and 
\be\|a_k\|_{\dot H^{1/2}}+\|a_k\|_{L^{\infty}}\le  C \|(-\Delta)^{1/4}u_k\|_{L^2(\R)}.\ee
   \end{Theorem}
From Theorem \ref{th-I.3} we will deduce that in a neck region (centered for simplicity of notation in $x=0$)  we have
\ \begin{equation}\label{limitneck2}
\lim_{\Lambda\to \infty}\lim_{k\to \infty}\left(\int_{B(0,\frac{R_k}{2\Lambda_k})\setminus B(0,\Lambda_k r_r)}|((-\Delta)^{1/4}u_k)^{-}|^2 dx\right)^{1/2}=0.
\end{equation}
where for a given  function $f\colon\R^n\to\R^m$ we denote by
$f^+$ and $f^-$ respectively the symmetric and antisymmetric part of $f$ with respect to the origin, i.e.
$f^+(x)=\frac{f(x)+f(-x)}{2}$ and $f^-(x)=\frac{f(x)-f(-x)}{2}$. It remains to find a link between the symmetric and the antisymmetric part of $(-\Delta)^{1/4}u_k$.
To this purpose we  
  make use of the {\bf Pohozaev type formula} theorem~\ref{th-I.4}, which is one of the main result of this paper.

We explain below in some steps how formula \rec{I-identityPR} permits to get an information of the $L^2$ norm of $((-\Delta)^{1/4}u_k)^-$.\par
 
{$\bullet$} We observe that we can rewrite the l.h.s and r.h.s  of \rec{I-identityPR} respectively as follows
\begin{eqnarray*} 
 \left| \int_{x\in\R} \frac{x^2-t^2}{(x^2+t^2)^2}u^+(x)dx\right|^2
  &=&t^{-2}\left\vert \int_{x\in\R} (-\Delta)^{-1/4}\left[\frac{x^2-1}{(x^2+1)^2}\right]((-\Delta)^{1/4}u)^{+}(xt)dx\right|^2.
  \end{eqnarray*}
  \begin{eqnarray*}\label{I-identityPR3}
\left| \int_{x\in\R}  \frac{2xt}{(x^2+t^2)^2}  u(x)dx\right|^2
  &=&t^{-2}\left\vert \int_{x\in\R} (-\Delta)^{-1/4}\left[\frac{2x}{(x^2+1)^2}\right]((-\Delta)^{1/4}u)^{-}(xt)dx\right|^2.
  \end{eqnarray*}
  Therefore we can rewrite the formula \rec{I-identityPR} as
  \begin{eqnarray}\label{I-identityPR4}
&& \left| \int_{x\in\R} (-\Delta)^{-1/4}\left[\frac{x^2-1}{(x^2+1)^2}\right]((-\Delta)^{1/4}u)^{+}(xt)\right\vert^2\\
&&~~~~~=\left|\int_{x\in\R} (-\Delta)^{-1/4}\left[\frac{2x}{(x^2+1)^2}\right]((-\Delta)^{1/4}u)^{-}(xt)dx\right|^2.\nonumber
  \end{eqnarray}\par
{$\bullet$} We prove that
$$M^{+}[v](t):=\int_{x\in\R} (-\Delta)^{-1/4}\left[\frac{x^2-1}{(x^2+1)^2}\right]v(xt)dx$$ and $$M^{-}[v](t):=\int_{x\in\R} (-\Delta)^{-1/4}\left[\frac{2x}{(x^2+1)^2}\right]v(xt)dx$$  are
isomorphisms from $L^{p}_+(\R)$ onto $L^{p}_+(\R)$ (resp . $L^{p}_-(\R)$ onto $L^{p}_-(\R)$) for every $p>1$ and by interpolation from $L^{2,1}_+(\R)$ onto  $L^{2,1}_+(\R)$(resp. $L^{2,1}_-(\R)$ onto  $L^{2,1}_-(\R)$). Here $L^p_+$ is the space of even $L^p$ functions and $L^p_-$ the space of odd $L^p$ functions. \par
 {$\bullet$}    By 
plugging  the development \rec{devel}  into \rec {I-identityPR4} we  obtain an information also of the asymptotic behaviour of the $L^2$ norm of $((-\Delta)^{1/4}u_k)^+$  in
an annular region, namely for all $\Lambda>1$ with  $\Lambda r_k<(2\Lambda)^{-1}R_k$  and $x\in B(0,(2\Lambda)^{-1}R_k)\setminus B(0,\Lambda r_k)$ 
\be\label{devel2}
\int_{B\left(0,\frac{R_k}{2\Lambda_k}\right)\setminus B(0,\Lambda_k r_r)}|((-\Delta)^{1/4}u_k)^{-}|^2 dx=\tilde h_k+\tilde g_k,\ee
where $\tilde g_k\in L^2(\R)$ with
$\displaystyle \limsup_{k\to\infty}\|\tilde g_k\|_{L^2} =0$ and $$\limsup_{\Lambda\to\infty}\limsup_{k\to+\infty} \|\tilde h_k\|_{L^{2,1}(A_{\Lambda r_k,(2\Lambda)^{-1}R_k})}<\infty.$$
{$\bullet$}  By combining \rec{devel} and \rec{devel2} we get   
\begin{Theorem}\label{th-I.5}
Under the assumptions of Theorem \ref{th-I.3}, we have 
\be
\limsup_{\Lambda\to\infty}\limsup_{k\to+\infty}\|(-\Delta)^{1/4}u_k\|_{L^2(B(0,(2\Lambda)^{-1}R_k)\setminus B(0,\Lambda_k r_r))}=0.\ee
  \end{Theorem}
 
We would like to establish a link between the above results and the result obtained by the second and third author \cite{LR} in the framework of harmonic maps in $2D$ with values into a closed manifold and that can easily be extended to
horizontal harmonic maps. 
Also in this case, the strategy has been to prove an  $L^{2,1}$-estimate on the angular part of the gradient (which play the role of the antisymmetric part $(-\Delta)^{1/4}u$). Precisely
suppose we have  a  sequence  $u_k\in{\frak{H}}^{1} (D^2)$  of horizontal harmonic maps such that $\|\nabla u_k\|_{L^2(D^2)}\le C.$ Then they satisfy
\be\label{I-Schr}-\Delta u_k=\Omega_k(P_T)\nabla u~~\mbox{in ${\cal{D}}^{\prime}(D^2)$},\ee
where  $\Omega\in L^{2}(D^2,{\R}^2\otimes so(m))$ is an antisymmetric potential. In \cite{LR} the authors found the following asymptotic development for $\nabla u_k$ in a annular domain $B(0,R_k)\setminus B(0,r_k)$ where 
$\|\Omega_k\|_{L^2(B(0,R_k)\setminus B(0,r_k)}$ is {\em small}:
\begin{equation}\label{develharm}
   \nabla u_k(x)= \frac{a_k(|x|)\overrightarrow{c_{r_k}}}{|x|^{1/2}} +  h_k,
    \end{equation}
where  $a_k$ is radial, $a_k\in L^{\infty}(B(0,R_k)\setminus B(0,r_k)),$ $h_k\in L^{2,1}(A_{\Lambda r_k,(2\Lambda)^{-1}R_k})$ (for every $\Lambda>1$ such that $\Lambda r_k<(2\Lambda)^{-1}R_k$),
\be\limsup_{\Lambda\to\infty}\limsup_{k\to+\infty} \|h_k\|_{L^{2,1}(A_{\Lambda r_k,(2\Lambda)^{-1}R_k})}<\infty,~~\|a_k\|_{L^{\infty}}\le C \|\nabla u_k\|_{L^2}.\ee
 From \rec{develharm} they deduce that 
 \begin{equation}\label{l21harm}
\limsup_{\Lambda\to\infty}\limsup_{k\to+\infty}  \left\Vert \frac{1}{\rho}\frac{\partial u_k}{\partial\theta}\right\Vert_{L^{2,1}(A_{\Lambda r_k,(2\Lambda)^{-1}R_k})}<\infty.
\end{equation}
 They also prove that the $L^{2,\infty}$ of  $\nabla u_k$ is arbitrary small in degenerating annuli, namely
\be\label{linftyharm}
 \limsup_{\Lambda\to\infty}\limsup_{k\to+\infty}  \|\nabla u_k\|_{L^{2,\infty}(A_{\Lambda r_k,(2\Lambda)^{-1}R_k})}=0.\ee
 By combining \rec{l21harm} and \rec{linftyharm} they obtain
 \begin{equation}
\limsup_{\Lambda\to\infty}\limsup_{k\to+\infty}  \left\Vert \frac{1}{\rho}\frac{\partial u_k}{\partial\theta}\right\Vert_{L^{2}(A_{\Lambda r_k,(2\Lambda)^{-1}R_k})}=0.
\end{equation}
The Pohozaev identity in $2D$ (\ref{PI}) obtained from \rec{IharmequationR2} implies that there is no loss of energy in the neck region:
 \begin{equation}
\limsup_{\Lambda\to\infty}\limsup_{k\to+\infty}  \|\nabla u_k\|_{L^{2}(A_{\Lambda r_k,(2\Lambda)^{-1}R_k})}=0.
\end{equation}

In a last part, we give a counter example to quantization for general sequence of Schr\"odinger equation with antisymmetric potential. For the Laplacian it has been donne in \cite{LR}. Hence here we gives an example of sequence of $u_k$ satisfying
$$(-\Delta)^{\frac{1}{2} }u_k =\Omega_k u_k +{\Omega_1}_k u_k$$
with $\Omega_k$ antisymmetric whose $L^2$-norm goes to zero in a neck region, and the $L^{2,1}$-norm  ${\Omega_1}_k$  goes to zero in the neck region, despite the $L^2$-norm of $u_k$ remains bounded from below. As in the case of  the Laplacian in \cite{LR}, we take a modification of the Green function. Since our functions will be even, this insures that that no Pohozaev identity can occur.
  
\subsection{Definitions and Notations} \label{secdef}
In this section we introduce some definitions and notations  that will be used in the paper.\par
 \begin{enumerate}

\item Given $\varphi\colon\R\to\C$ we denote by $\hat \varphi$ or ${\cal{F}}[\varphi]$ the Fourier transform of $\varphi$ and
by $\check\varphi$ or ${\cal{F}}^{-1}[\varphi]$ the inverse Fourier transform of $\varphi.$
\item Given $u,v\in\R^n$ we denote by $\langle u, v\rangle $ the standard scalar product of $u$ and $v$.
\item Given $z\in \C$ we denote by $\Re{z}$ and $\Im{z}$ the real and imaginary part of $z$.
\end{enumerate}
Next we introduce some functional spaces.\\

First, we will note $L^p_+$  the space of even $L^p$ functions and $L^p_-$ the space of odd $L^p$ functions.\\

We denote by $L^{2,\infty}(\R^n)$ the space of measurable functions $f$ such that
$$
\sup_{\lambda>0}\lambda |\{x\in\R^n~: |f(x)|\ge\lambda\}|^{1/2}<+\infty\,,
$$
and $L^{2,1}(\R^n)$ is the  space of measurable functions satisfying
$$\int_{0}^{+\infty}|\{x\in\R~: |f(x)|\ge\lambda\}|^{1/2} d\lambda<+\infty \,\,.$$
We can check that $L^{2,\infty}$ and  $L^{2,1}$ forms a duality pair. \\

We denote by ${\mathcal H}^1({\R}^n)$   the Hardy space which is the space of $L^1$
functions $f$ on ${\R}^n$satisfying 
\[
\int_{{\R}^n}\sup_{t>0}|\phi_t\ast f|(x)\ dx<+\infty\quad ,
\]
where $\phi_t(x):=t^{-n}\ \phi(t^{-1}x)$ and   $\phi$ is some function in the Schwartz space ${\mathcal S}({\R}^n)$ satisfying $\int_{{\R}^n}\phi(x)\ dx=1$.\\

Finally for every $s\in\R$ and $q>1$ we denote  by 
$\dot W^{s,q}(\R^n)$ the fractional Sobolev space
 $$\{f\in{\cal{S}}^{\prime}:~{\cal{F}}^{-1}[|\xi|^s|{\cal{F}}[f]]\in L^q(\R^n)\}\,.$$
For more properties on the Lorentz spaces, Hardy space ${\mathcal H}^1$ and   fractional Sobolev spaces we refer to \cite{Gra1} and \cite{Gra2}.

 The paper is organized as follows.  
 In Section \ref{nonlocalSchr} we analyze the asymptotic behaviour of the solutions of special nonlocal Schr\"odinger systems with a $L^2$ antisymmetric potential in degenerate annular domain.
 In Section \ref{PoId} we describe two new Pohozaev type formulae, one on $\R$ and one in $S^1$, which can be obtained one from the other by means of the stereographic projections and   which are satisfied in particular by $1/2$-harmonic maps respectively on $\R$ and on $S^1.$ We also compare such formulae with a Pohozev formula in $2D$ which is satisfied by harmoinc maps. In Section 4 we deduce from the results of Section 2 and 3 the quantization of horizontal $1/2$-harmonic maps. Finally in Section 5 we describe an example showing that in general we cannot have quantization for solutions of a general nonlocal Schr\"odinger systems with a $L^2$ antisymmetric potential.

\section{Preliminary results on nonlocal Schr\"odinger system}\label{nonlocalSchr}
 
  As we will see later horizontal $1/2$-harmonic maps in $1$-D satisfy special  
 nonlocal system of the form
\begin{equation}\label{modelsystem}
(-\Delta)^{1/4} v=\Omega v+\Omega_1 v+{\cal{Z}}(Q,v)+g \end{equation}
where $v\in L^2(\R)$, $Q\in  H^{1/2}(\R)$,  $\Omega\in L^2(\R,so(m))$, $\Omega_1\in L^{2,1}(\R)$, $g\in L^1(\R)$ and ${\cal{Z}}\colon  H^{1/2}(\R)\times L^2(\R)\to {\cal{H}}^1(\R)$  satisfies the following stability  property:    if $Q,P\in{ \dot {H}}^{1/2}(\R)\cap L^{\infty}(\R)$, $v\in L^{2}$ then  
\begin{equation}\label{decZ}
P{\cal{Z}}(Q,v)={{{\cal{A}}_{{\cal{Z}}}}}(P,Q)v+J_{{\cal{Z}}}(P,Q,v),
\end{equation}
where
$$
 \|{\cal{A}}_{{\cal{Z}}}(P,Q)\|_{L^{2,1}(\R)}\leq C\|(-\Delta)^{1/4}[P] \|_{L^2(\R)}\|(-\Delta)^{1/4}[Q] \|_{L^2(\R)},$$
 and
 $$
\|{{J}}_{{\cal{Z}}}(P,Q,v)\|_{{\cal{H}}^1(\R)}\le C \left(\|(-\Delta)^{1/4}[P] \|_{L^2(\R)}+\|(-\Delta)^{1/4}[Q] \|_{L^2(\R) }\right)\|v\|_{L^2(\R)}.$$
Actually ${\cal{Z}}$ is a linear combination of the following {\em pseudo-differential operators}:
\begin{equation}\label{opT}
T(Q,v):=(- \Delta)^{1/4}(Qv)-Q(- \Delta)^{1/4} v+ (- \Delta)^{1/4} Q v\,\end{equation}
and
\begin{equation}
\label{opS}
S(Q,v):=(-\Delta)^{1/4}[Qv]-{\cal{R}} (Q{\cal{R}}(-\Delta)^{1/4} v)+{\cal{R}}((-\Delta)^{1/4} Q{\cal{R}} v )
\end{equation}
 \begin{equation}\label{opF}
 F(Q,v):={\mathcal{R}}[Q]{\mathcal{R}}[v]-Qv.
 \end{equation}
 \begin{equation}\label{opOm}
 \Lambda(Q,v):=Qv+{\mathcal{R}}[Q{\mathcal{R}}[v]]\,.
 \end{equation}
 In  \cite{DL2} the first author impove the estimates on the operators $T,S$ obtained in \cite{DLR1}:
 \begin{Theorem}\label{comm2}
 Let $v\in L^2(\R),$ $Q\in \dot{H}^{1/2}( \R)$. Then $T(Q,v),S(Q,v)\in {\cal{H}}^{1}( \R)$  and
\begin{equation}\label{commest2}
\|T(Q,v)\|_{{\cal{H}}^{1}( \R)}\le C\|Q\|_{\dot{H}^{1/2}( \R)}\|v\|_{L^2( \R)}\,.\end{equation}
\begin{equation}\label{commest3} 
\|S(Q,v)\|_{{\cal{H}}^{1}( \R)}\le C\|Q\|_{\dot{H}^{1/2}( \R)}\|v\|_{L^2( \R)}\,.~~~\hfill \Box
\end{equation}
\end{Theorem}
\medskip
As a consequence of the Coifman-Rochberg-Weiss   estimate \cite{CRW}  we also have
\begin{equation}\label{estFhardy}
\|F(f,v)\|_{  {\cal{H}}^1(\R)}\le C \|f\|_{L^2(\R)}\|v\|_{L^{2}(\R)}\,.
\end{equation}
for $f,v\in L^2$.\par
In \cite{DLR3} the first and the third author show that the stability property \rec{decZ} holds for the operators \rec{opT}, \rec{opS},\rec{opF}.

Moreover the authors have shown  that if the $L^2$ norm of $\Omega$ is ``small" then   the system \rec{modelsystem} is equivalent to
a conservation law:
 
 \begin{Theorem}[Theorem 3.11 in \cite{DLR3}]\label{conservation}
 Let $v\in L^{2}(\R,\R^m)$ be a solution of \rec{modelsystem},
 where $\Omega\in L^2(\R,so(m))$, $\Omega_1\in L^{2,1}(\R)$,   ${\cal{Z}}(Q,v)\in {\cal{H}}^1$ for every $Q\in{H^{1/2}},$ $v\in L^2$ with
\begin{eqnarray*}
\|{\cal{Z}}(Q,v)\|_{{\cal{H}}^1}&\le& C \| Q\|_{{H^{1/2}}}\|v\|_{L^2},
\end{eqnarray*}
and ${\cal{Z}}(Q,v)$ satisfies \rec{decZ}. There exists $\gamma>0$ such that if $\|\Omega\|_{L^2}<\gamma$, then 
 there exist    $A=A(\Omega,\Omega_1,Q)\in L^{\infty}\cap H^{1/2}(\R,M_m({\R})))$ and an operator $ B=B(\Omega,\Omega_1,Q)\in  H^{1/2}(\R)$ such that
\begin{eqnarray}
\|A\|_{ H^{1/2}}+\|B\|_{ H^{1/2}}&\le& C(\|\Omega\|_{L^2}+\|\Omega_1\|_{L^{2,1}} +\| Q\|_{{H^{1/2}}})\\
dist(\{A,A^{-1} \},SO(m))&\le& C(\|\Omega\|_{L^2}+\|\Omega_1\|_{L^{2,1}} +\| Q\|_{{H^{1/2}}})
\end{eqnarray}
and 
\begin{equation}\label{conslaw}
(-\Delta)^{1/4}[Av]={\cal{J}}(B,v)+Ag,
\end{equation}
where ${\cal{J}}$ is a linear operator in $B,v$, ${\cal{J}}(B,v)\in {\cal{H}}^1(\R)$ and 
\begin{equation}\label{J}
\|{\cal{J}}(B,v)\|_{{\cal{H}}^1(\R)}\le C \| B \|_{H^{1/2}}\|v\|_{L^2}\,.
\end{equation}
 \end{Theorem}
 
In the following result we show that  if the $L^2$-norm of the antisymmetric potential $\Omega$ in \rec{modelsystem} is {\em small} in an annular domain 
$A_{r,R}:=B(0,R)\setminus B(0,r)$ then a precise asymptotic development inside $A_{r,R}$  holds for the solution $v$ to \rec{modelsystem}.\par
We assume for simplicity that $g\equiv 0.$
 
 \begin{Proposition}\label{decomposition}
 Let $v\in L^{2}(\R,\R^m)$ be a solution of 
\begin{equation}\label{modelsystem2}
(-\Delta)^{1/4} v=\Omega v+\Omega_1 v+{\cal{Z}}(Q,v) \end{equation}
where $\Omega\in L^2(\R,so(m))$, with $\|\Omega\|_{L^2(B(0,R)\setminus B(0,r))}<\gamma$, $\Omega_1\in L^{2,1}(\R)$,  
${\cal{Z}}(Q,v)\in {\cal{H}}^1$ for every $Q\in{H^{1/2}},$ $v\in L^2$ with
\begin{eqnarray*}
\|{\cal{Z}}(Q,v))\|_{{\cal{H}}^1}&\le& C \| Q\|_{{H^{1/2}}}\|v\|_{L^2},
\end{eqnarray*}
and ${\cal{Z}}(Q,v)$ satisfies \rec{decZ}
 Then 
 there exists $A=A(\Omega,\Omega_1,Q)\in H^{1/2}(\R,M_m({\R})))$ with  
 \begin{eqnarray}
\|A\|_{ H^{1/2}}  +\|A^{-1}\|_{ H^{1/2}} &\le& C(\|\Omega\|_{L^2}+\|\Omega_1\|_{L^{2,1}} +\| Q\|_{{H^{1/2}}})\\
dist(\{A,A^{-1} \},SO(m))&\le& C(\|\Omega\|_{L^2}+\|\Omega_1\|_{L^{2,1}} +\| Q\|_{{H^{1/2}}}),\nonumber
\end{eqnarray}
 and $\overrightarrow{c_{r}}\in\R^m$ such that for every $x\in B(0,R)\setminus B(0,r)$ and for every $\Lambda>2$ with $\Lambda r< (2\Lambda)^{-1}R$ it holds
\begin{equation}\label{decompv}
   v(x)=(A^{-1}(x))^{+}\overrightarrow{c_{r}}\frac{1}{|x|^{1/2}}+h(x)+g(x)\end{equation}
 with $h\in L^{2,1}(A_{\Lambda r,(2\Lambda)^{-1}R})$, 
 $\|h\|_{L^{2,1}(A_{\Lambda r, (2\Lambda)^{-1}R})}\le   C\|v\|_{L^2}$ where $C$ is a positive constant independent of $r,R,\Lambda$ and $g\in L^2(\R)$ with
 $\|g\|_{L^2(\R)}\le C |\overrightarrow{c_{r}}|.$
   \end{Proposition}
   
{\bf Proof of Proposition \ref{decomposition}.}
 We split $\Omega$ as follows:
$$
\Omega=\Omega^{r}+\Omega^{R}+  \Omega^{R,r}$$
where 
\begin{eqnarray*}
  \Omega^{R,r} &=&\11_{B(0,R)\setminus B(0,r)}\Omega\\
\Omega^{r}&=&\11_{ B(0,r)}\Omega\\
  \Omega^{R}&=&\11_{ B^c(0,R)}\Omega.
  \end{eqnarray*}
  We write the system \rec{modelsystem2} in the following form
  \begin{eqnarray}\label{modelsystem3}
(-\Delta)^{1/4} v&=&\Omega^{R,r} v+\Omega_1 v+{\cal{Z}}(Q,v)+h_{r}+h_R\end{eqnarray}
where
$$
h_{r}=\Omega^{r}\, v,~~\mbox{and} ~~h_{R}=\Omega^{R}\, v .$$
 
Observe that $ h_{r}, h_{R}\in L^1(\R)$ with $supp(h_R)\subseteq B^c(0,R),$ $supp(h_{r})\subseteq B(0,r)$, and $\|h_{r}\|_{L^1}+\|h_{R}\|_{L^1}\le C(\|v\|_{L^2}).$
From Theorem \ref{conservation} there exist    $A=A(\Omega,\Omega_1,Q)\in \dot H^{1/2}(\R,M_m({\R}))$ and an operator $ B=B(\Omega,\Omega_1,Q)\in  \dot H^{1/2}(\R)$ such that
\begin{eqnarray}
\|A\|_{ \dot H^{1/2}} +\|A^{-1}\|_{ \dot H^{1/2}}+\|B\|_{ \dot H^{1/2}}&\le& C(\|\Omega\|_{L^2}+\|\Omega_1\|_{L^{2,1}} +\| Q\|_{{\dot H^{1/2}}})\\
dist(\{A,A^{-1} \}, SO(m))&\le& C(\|\Omega\|_{L^2}+\|\Omega_1\|_{L^{2,1}} +\| Q\|_{{\dot H^{1/2}}})
\end{eqnarray}
and 
\begin{equation}\label{conslaw2}
(-\Delta)^{1/4}[Av]={\cal{J}}(B,v)+Ah_{r}+Ah_{R},
\end{equation}
with ${\cal{J}}(B,v)$ satisfying \rec{J}.
We write
$v=v_1+v_2+v_3+v_4$ where
\begin{eqnarray} 
(-\Delta)^{1/4}[Av_1]&=&{\cal{J}}(B,v);\label{v1}\\
(-\Delta)^{1/4}[Av_2]&=& Ah_{r}-\left(\int_{B(0,r)}Ah_{r} dx\right)\delta_0;\label{v2}\\
(-\Delta)^{1/4}[Av_3]&= &Ah_{R},\label{v3}\\
(-\Delta)^{1/4}[Av_4]&=&  \left(\int_{B(0,r)}Ah_{r} dx\right)\delta_0;\label{v4}.
 \end{eqnarray}
 Let $\Lambda>2$ be such that $\Lambda r< (2\Lambda)^{-1}R.$
 We are going to estimate the $L^{2,1}(B(0,(2\Lambda)^{-1}R)\setminus B(0, \Lambda r))$ norm of $v_1,v_2,v_3$.\par
 {\bf 1. Estimate of $v_1$\,.}\par 
 By the properties of $\cal{J}$ we have
 \begin{eqnarray}\label{estv1}
 \|Av_1\|_{L^{2,1}(\R)} &\leq& C\|(-\Delta)^{-1/4}[J(B,v)]\|_{L^{2,1}(\R)}\\
&\leq& C \| B \|_{\dot H^{1/2}}\|v\|_{L^2}.
\end{eqnarray}

 {\bf 2. Estimate of $v_2$\,.}\par
 We set $\phi_{r}:=Ah_{r}-\left(\int_{B(0,r)}Ah_{r} dx\right)\delta_0$ and take $x\in B(0,(2\Lambda)^{-1}R)\setminus B(0, \Lambda r)\,.$
By using the fact that $\int_{B(0,r)} \phi_{r}(y) dy=0$, we get  
 \begin{eqnarray*}
Av_2(x)&=&\int_{B(0,r)}\phi_{r}(y)\frac{1}{|x-y|^{1/2}} dy= \frac{1}{|x|^{1/2}}\sum_{k=0}^{+\infty} \frac{c_{k}}{|x|^k}\int_{B(0,r)}\phi_{r}(y)y^k dy\\
 &=& \sum_{k=1}^{+\infty} \frac{c_k}{|x|^{k+1/2}}\int_{B(0,r)}\phi_{r}(y)y^k dy\,,
 \end{eqnarray*}
 where the $c_k=\frac{1\times 3\times\dots (2k-1)}{2\times 4\times \dots \times 2k}$.\\
 
 We observe that for every $k\ge 1$, $\frac{1}{|x|^{k+1/2}}\in L^{2,1}(B(0,(2\Lambda)^{-1}R)\setminus B(0, \Lambda r))$ for all $\Lambda>1$ such that $\Lambda r< (2\Lambda)^{-1}R$, with
 $$\left\|\frac{1}{|x|^{k+1/2}}\right\|_{ L^{2,1}( B(0,(2\Lambda)^{-1}R)\setminus B(0, \Lambda r))}\precsim \left(\frac{1}{ \Lambda r}\right)^{k}\,.$$
 Therefore
 \begin{eqnarray*}
    \|v_2\|_{ L^{2,1}( B(0,(2\Lambda)^{-1}R)\setminus B(0,  \Lambda r))}&\leq&C \|\phi_{r}\|_{L^1}
 \sum_{k=1}^{+\infty} \frac{c_k}{( \Lambda r)^k}r^k\\
 & \leq &C \|\phi_{r}\|_{L^1}
 \sum_{k=1}^{+\infty} \frac{c_k}{\Lambda ^k}<+\infty\\
 &\leq& C \|\Omega\|_{L^{2}}\|v\|_{L^{2}}\left( \sum_{k=1}^{+\infty} \frac{c_k}{\Lambda ^k}<+\infty\right).
  \end{eqnarray*}
  {\bf 3. Estimate of $v_3$\,.}\par
 Since $Ah_{R}$ has support in $B^c(0,R)$ then by interpolation one gets that 
$A v_3\in L^{2,1}(B(0,(2\Lambda)^{-1}R)\setminus B(0, \Lambda r)))$ and
$$\|A v_3\|_{ L^{2,1}(B(0,(2\Lambda)^{-1}R)\setminus B(0, \Lambda r))}\le C \|h_R\|_{L^1}\le C\|\Omega\|_{L^{2}}\|v\|_{L^{2}})\,.$$
{\bf 3. Estimate of $v_4$\,.}\par
We have
$$v_4(x)=A^{-1}(x)\overrightarrow{c_{r}}\frac{1}{|x|^{1/2}},$$
where
$$\overrightarrow{c_{r}}=\int_{B(0,r)}Ah_{r} dx.$$
Since $A^{-1} \in \dot H^{1/2}(\R)$, it verifies the following estimate (see lemma 33.1 of \cite{Tartar}):
   \begin{eqnarray}\label{estA}
   \left\|\frac{A^{-1}(x)-A^{-1}(-x)}{|x|^{1/2}}\right\|_{L^2(\R)}&\le& C\|A^{-1}\|_{\dot H^{1/2}}\end{eqnarray}
   Therefore 
by combining the $L^{2,1}$ estimates on $ B(0,(2\Lambda)^{-1}R)\setminus B(0, \Lambda r)$ of  $v_1,v_2,v_3$ we obtain
\begin{eqnarray}\label{estl21v}
 v(x)&=&(A^{-1}(x))^{+}\overrightarrow{c_{r}}\frac{1}{|x|^{1/2}}  +(A^{-1}(x))^{-}\overrightarrow{c_{r}}\frac{1}{|x|^{1/2}} +h(x)
 \end{eqnarray}
with $h(x)=v_1+v_2+v_3$ and $\|h\|_{L^{2,1}(B(0,(2\Lambda)^{-1}R)\setminus B(0, \Lambda r))}\le C\|v\|_{L^2(\R)},$ $\overrightarrow{c_{r}}=\int_{B(0,r)}A(\Omega^{r}\, v)dx.$
Finally from \rec{estA} it follows that $g(x)={\rm asymm}{(A^{-1}(x))}\overrightarrow{c_{r}}\frac{1}{|x|^{1/2}}$ satisfies
$g\in L^2(\R)$ with $\|g\|_{L^2(\R)}\leq C |\overrightarrow{c_{r}}|.$
We can conclude the proof.\hfill$\Box$\par

\begin{Lemma}[$\varepsilon$-regularity ]\label{epsreg}
Let $v\in L^2(\R,\R^m)$ be a solution of
 \begin{equation}
(-\Delta)^{1/4} v=\Omega v+\Omega^1   v+{\cal{Z}}(Q ,v ).\end{equation}
Then there exists $\varepsilon_0>0$ such that if 
 \begin{equation}\label{epsregcond}\sum_{j\ge 0} 2^{-j/2}(\|\Omega\|_{L^2(B(x,2^{j}r)) } +\|\Omega^1\|_{L^{2,1}(B(x,2^{j}r)}+\|(-\Delta)^{1/4}Q\|_{L^2(B(x,2^{j}r)})\le \varepsilon_0\,,\end{equation}
then there is $p>2$ ( independent of $u$) such that for every   $y\in B(x,r/2)$ we have
\begin{equation}
\left(r^{p/2-1}\int_{B(y,r/2)}|v|^p dx\right)^{1/p}\le C\sum_{j\ge 0} 2^{-j/2}\left(\|\Omega\|_{L^2(B(x,2^{j}r)) } +\|\Omega^1\|_{L^{2,1}(B(x,2^{j}r)}+\|(-\Delta)^{1/4}Q\|_{L^2(B(x,2^{j}r)}\right),
\end{equation}
where $C>0 $ depends on $\|v\|_{L^2}$\,.
\end{Lemma}

  Finally we have
  
\begin{Lemma}[$L^{2,\infty}$ estimates]\label{neck2}
   There exists  $ \delta >0$ such that for any  $ v\in L^2(\R)$ solution of \rec{modelsystem} and $\lambda,\Lambda>0$  with $\lambda<(2\Lambda)^{-1} $ satisfying
 $$ \sup_{\rho\in[\lambda,{(2\Lambda)}^{-1}]} \left[\|\Omega\|_{L^{2}(B(0,2\rho)\setminus B(0,\rho))}+\|(-\Delta)^{1/4})Q\|_{L^{2}(B(0,2\rho)\setminus B(0,\rho))}+\|\Omega^1\|_{L^{2,1}(B(0,2\rho)\setminus B(0,\rho))}\right]\le \delta $$
then
\begin{eqnarray}\label{linftyest}
&& \| v\|_{L^{2,\infty}(B(0,(2\Lambda)^{-1})\setminus B(0,\lambda )}\le  C\sup_{\rho\in[\lambda,{(2\Lambda)}^{-1}]} \left[\|\Omega\|_{L^{2}(B(0,2\rho)\setminus B(0,\rho))}\right.\\
&&\left.+\|(-\Delta)^{1/4}Q\|_{L^{2}(B(0,2\rho)\setminus B(0,\rho))}+\|\Omega^1\|_{L^{2,1}(B(0,2\rho)\setminus B(0,\rho))}\right] \,.\nonumber
\end{eqnarray}
 where $C$ is  independent of $\lambda,\Lambda\,.$
 \end{Lemma}
  The proof of Lemmae \ref{epsreg}  and  \ref{neck2}  are the same of Lemma 2.1  and Lemma 3.2 in \cite{DL2} and therefore we omit them.
 
We  show now a point removability type result for solution of \rec{modelsystem}.
  \begin{Proposition}\label{removsing}[ Point removability]
  Let $v \in L^2(\R)$ be a solution of \rec{modelsystem} in  ${\cal{D}}^{\prime} (\R\setminus\{a_1,\ldots,a_{\ell}\})$. Then it is a solution
  of \rec{modelsystem} ${\cal{D}}^{\prime} (\R)$.
      \end{Proposition}
   {\bf Proof of Proposition \ref{removsing}\,.}
   The fact that
    $$
(-\Delta)^{1/4} v=\Omega v+\Omega^1   v+{\cal{Z}}(Q ,v ),
 \quad \mbox{in ${\cal{D}}^{\prime} (\R\setminus\{a_1,\ldots,a_{\ell}\})$}$$
implies that
the distribution
 $\phi:= (-\Delta)^{1/4} v-\Omega v-\Omega^1   v-{\cal{Z}}(Q ,v )$  is of  order $p=1$ and supported in $\{a_1,\ldots,a_{\ell}\}$.
 Therefore by Schwartz Theorem (see 6.1.5.of \cite{Bony}) one has
 $$\phi =\sum_{|\alpha|\le 1} c_{\alpha}\partial^{\alpha} \delta_{a_i}\,.
 $$
 Since $\phi\in \dot L^1(\R) $,      then  the above implies that $c_{\alpha}=0$ and thus
  $v$ is a solution of \rec{modelsystem} in {  ${\cal{D}}^{\prime} (\R )$}\,. 
 We conclude the proof of Proposition  \ref{removsing}\,.~\hfill$\Box$
  
    \bigskip
    \section{Pohozaev Identities}\label{PoId}
    In this section we provide some new Pohozaev identities in $\R$ and $S^1$ which are in particular satisfied by horizontal $1/2$ harmonic in ${\frak{H}}^{1/2} ({\R})$.
       
  \subsection{Pohozaev Identities for $(-\Delta)^{1/2}$ in $\R$}\label{Poh}
We first determine the solution of
\begin{equation}\label{solfund1}
\left\{\begin{array}{cc}
\partial_t G+(-\Delta)^{1/2} G=0 & t>0\\
G(0,x)=\delta & t=0\,.
\end{array}\right.
\end{equation}
In order to solve \rec{solfund1} we consider the Fourier Transform of $G$ with respect to $x$
$$\hat G(t,\xi)=\int_{\R} G(t,x) e^{-ix\xi} dx$$
The function $\hat G$ satisfies
the problem
\begin{equation}\label{solfund2}
\left\{\begin{array}{cc}
\partial_t \hat G+|\xi|\hat G=0 & x\in\R,\, t>0\\
\hat G(0,x)=\hat\delta=1 & x\in\R,\,t=0\,.
\end{array}\right.
\end{equation}
The solution of \rec{solfund2} is
$\hat G(t,\xi)=e^{-t|\xi|}$.
Then $$G(t,x)=\frac{1}{2\pi}\int_{\R} \hat G(t,\xi)e^{ix\xi} d\xi=\frac{1}{\pi} \frac{t}{x^2+t^2}.$$
The following equalities hold
$$\partial_t G=\frac{1}{\pi}\frac{x^2-t^2}{(t^2+x^2)^2},~~\partial_x G=-\frac{1}{\pi}\frac{2xt}{(t^2+x^2)^2}.$$
 
\begin{Theorem}{[Case on $\R$]}\label{Poho1}
Let $u\in \dot H^{1/2}_{loc}(\R,\R^m)$ such that
\begin{equation}\label{harmequation1}
\left\langle \frac{du}{dx}, (-\Delta)^{1/2} u\right\rangle=0~~\mbox{a.e in $\R$.}
\end{equation}
Assume that
\begin{equation}\label{condR}
\int_{\R}|u-u_0|dx <+\infty, ~~\int_{\R}\left|\frac{du}{dx}(x)\right|dx<+\infty..
\end{equation}
Then the following identity holds
\begin{equation}\label{idpohozaevR}
 \left|\int_{x\in\R} \frac{x^2-t^2}{(x^2+t^2)^2}u(x)dx\right|^2=\left| \int_{x\in\R}  \frac{2xt}{(x^2+t^2)^2}  u(x)dx\right|^2.
\end{equation}
\end{Theorem}
{\bf Proof.}
We assume for simplicity that $u_0=0$.
Set $G_t(x):=G(t,x),$ multiply the equation \rec{harmequation1} by $xG_t(x)$ and we integrate. For every $j\in\{1,\ldots,m\}$ we get by Plancherel Theorem
\begin{eqnarray}\label{est1}
0&=&\int_{\R} x G_t  \frac{du^j}{dx} \,\overline{(-\Delta)^{1/2} u^j} dx =\int_{\R}{\mathcal{F}}\left[xG_t  \frac{du^j}{dx}\right]\overline{{\mathcal{F}}[(-\Delta)^{1/2} u^j]}dx\nonumber 
\\&=&
\int_{\R}{\mathcal{F}}[G_t]\ast {\mathcal{F}}\left[x \frac{du^j}{dx}\right] \, |\xi| \overline{{\mathcal{F}}[u^j] }d\xi  \\ &
= &\int_{\R}{\mathcal{F}}[G_t]\ast i\frac{d}{d\xi}{\mathcal{F}}\left[ \frac{du^j}{dx}\right] \, |\xi| \overline{{\mathcal{F}}[u^j] }d\xi\nonumber \\ &
= &\int_{\R}{\mathcal{F}}[G_t]\ast i\frac{d}{d\xi}\left(\xi  {\mathcal{F}}[u^j]\right)\, |\xi|\overline{{\mathcal{F}}[u^j] }d\xi \nonumber\\ &
= &
-\iint_{\R^2}\hat {G_t}(\xi-\eta)\frac{\partial}{\partial \eta}(\eta \hat{u}^j(\eta))|\xi| 
\overline{\hat{u}^j}(\xi) d\xi d\eta \nonumber\\
&=& \iint_{\R^2}\frac{\partial}{\partial \eta}\left[\hat {G_t}(\xi-\eta)\right] \eta \, |\xi|\hat{u}^j(\eta))
\overline{\hat{u}^j}(\xi) d\xi d\eta\nonumber \\
&=&\iint_{\R^2} te^{-t|\xi-\eta|}\frac{\xi-\eta}{|\xi-\eta|}\eta 
  |\xi| \hat{u}^j(\eta))
\overline{\hat{u}^j}(\xi) d\xi d\eta \nonumber
\end{eqnarray}
We symmetrize \rec{est1}  
\begin{eqnarray}\label{sym}
0&=&
 \iint_{\R^2} te^{-t|-\xi+\eta|}\frac{-\xi+\eta}{|-\xi+\eta|} \xi 
  (|\eta|) \hat{u}^j(\xi))
\overline{\hat{u}^j}(\eta) d\xi d\eta 
\end{eqnarray}
Then we sum \rec{est1} and \rec{sym} and we take the real part:
\begin{eqnarray}\label{real}
0&=&
 \iint_{\R^2} te^{-t|\xi-\eta|}\frac{\xi-\eta}{|\xi-\eta|}[\eta |\xi|-\xi
  |\eta|]\Re\left[ \hat{u}^j(\xi)\overline{\hat{u}^j} 
 \right] d\xi d\eta\nonumber\\
 &&=\iint_{\R^2} te^{-t|\xi-\eta|}\frac{\xi-\eta}{|\xi-\eta|}[\eta |\xi|-\xi
  |\eta|]\left[ a^j(\xi)a^j(\eta)+b^j(\xi)b^j(\eta)
 \right] d\xi d\eta
\end{eqnarray}

where for every $j$ we set 
$\hat{u}^j(\eta))=a^j(\eta)+i b^j(\eta).$ \par
 
We write 
 
\begin{eqnarray*}
a^{j,\pm}(\eta)&=&\frac{a^j(\eta)\pm a^j(-\eta)}{2}\\
b^{j,\pm}(\eta)&=&\frac{b^j(\eta)\pm b^j(-\eta)}{2}.
\end{eqnarray*}
We observe that $a^{j,+}, b^{j,+}$ are even and $a^{j,-}, b^{j,-}$ are odd. Moreover  since $u^j$ is real we also have
$b^{j,+}=0, a^{j,-}=0$, $b^{j,-}=b, a^{j,-}=a$.
We can also write
\begin{eqnarray}\label{estaplus}
a^{j,+}(\xi)&=&\int_{\R}\frac{u^j(x)+u^j(-x)}{2}\cos(x\xi)\,dx-i\int_{\R}\frac{u^j(x)+u^j(-x)}{2}\sin(x\xi)\,dx\nonumber\\ &&={\mathcal{F}}(u^{j,+})(\xi);
\end{eqnarray}
\begin{eqnarray}\label{estbminus}
b^{j,-}(\xi)&=&i\left\{\int_{\R}\frac{u^j(x)-u^j(-x)}{2}\cos(x\xi)\,dx-i\int_{\R}\frac{u^j(x)-u^j(-x)}{2}\sin(x\xi)\,dx\right\}\nonumber\\ &&=i{\mathcal{F}}(u^{j,-})(\xi)
\end{eqnarray}

We set
$$
Q(\alpha,\beta):=\iint_{\R^2}\frac{\xi-\eta}{|\xi-\eta|}te^{-t(|\xi-\eta|)}[\eta \, |\xi|-\xi \, |\eta|]\alpha(\eta)\beta(\xi) d\xi\, d\eta.$$
 It  follows that
\begin{eqnarray}\label{q}
0&=&Q(a^{j,+} ,a^{j,+}  )+Q(b^{j,-} ,b^{j,-} ).
 \end{eqnarray}

{\bf Estimate of $Q(a^{j,+} ,a^{j,+}  )$}
 
\begin{eqnarray}\label{aplus}
Q(a^{j,+} ,a^{j,+} )&=&\iint_{\R^2}   \frac{\xi-\eta}{|\xi-\eta|}te^{-t(|\xi-\eta|)}[\eta \, |\xi|-\xi \, |\eta|] a^{j,+} (\eta)a^{j,+} (\xi)d\xi\,d\eta\\
&=&
2 \int_{\xi>0} \int_{\eta<0} \eta\xi  \frac{\xi-\eta}{|\xi-\eta|}te^{-t(|\xi-\eta|)} a^{j,+} (\eta)a^{j,+} (\xi)d\xi\,d\eta\nonumber\\
&-& 2\int_{\xi<0} \int_{\eta>0} \eta\xi  \frac{\xi-\eta}{|\xi-\eta|}te^{-t(|\xi-\eta|)}  a^{j,+} (\eta)a^{j,+} (\xi)d\xi\,d\eta\nonumber\\
&=& 4\int_{\xi>0} \int_{\eta<0} \eta\xi  te^{-t(|\xi-\eta|)}  a^{j,+}(\eta) a^{j,+} (\xi)d\xi\,d\eta\nonumber
\\
&=&-4t\left(\int_{\xi>0}  \xi  e^{-t\xi}  a^{j,+} (\xi)d\xi\right)^2\nonumber\\
&=&-4t \left(\int_{\xi>0}  \xi  e^{-t\xi}  {\mathcal{F}}(u^{j,+})(\xi)d\xi\right)^2=-4t\left(\frac{1}{2}\int_{x \in\R}| \xi | e^{-t|\xi |}  {\mathcal{F}}(u^{j,+})(\xi)d\xi\right)^2\nonumber\\
&=&- t \left(\int_{\xi\in\R} \mathcal{F}{[(-\Delta)^{1/2} G](\xi)}{\mathcal{F}}(u^{j,+})(\xi)d\xi\right)^2\nonumber\\
&=&- t \left(\int_{x\in\R} \partial_{t} G(x,t)  (u^{j,+})(x)dx\right)^2\nonumber\\
&=&t\frac{1}{\pi^2} \left(\int_{x\in\R} \frac{x^2-t^2}{(x^2+t^2)^2}  (u^{j,+})(x)dx\right)^2\nonumber
\end{eqnarray}

{\bf Estimate of $Q(b^{j,-} ,b^{j,-}  )$}
\begin{eqnarray}\label{bminus}
Q(b^{j,-} ,b^{j,-} )&=&\iint_{\R^2}   \frac{\xi-\eta}{|\xi-\eta|}te^{-t(|\xi-\eta|)}[\eta \, |\xi|-\xi \, |\eta|] b^{j,-} (\eta)b^{j,-} (\xi)d\xi\,d\eta\\
&=&
2 \int_{\xi>0} \int_{\eta<0} \eta\xi  \frac{\xi-\eta}{|\xi-\eta|}te^{-t(|\xi-\eta|)} b^{j,-} (\eta)b^{j,-} (\xi)d\xi\,d\eta\nonumber\\
&-& 2\int_{\xi<0} \int_{\eta>0} \eta\xi  \frac{\xi-\eta}{|\xi-\eta|}te^{-t(|\xi-\eta|)}  b^{j,-} (\eta)b^{j,-} (\xi)d\xi\,d\eta\nonumber\\
&=& 4\int_{\xi>0} \int_{\eta<0} \eta\xi  te^{-t(|\xi-\eta|)}  b^{j,-} (\eta)b^{j,-} (\xi)d\xi\,d\eta\nonumber
\\
&=&4t\left(\int_{\xi>0}  \xi  e^{-t\xi}  b^{j,-} (\xi)d\xi\right)^2\nonumber\\
&=&4t \left(\int_{\xi>0}  \xi  e^{-t\xi} i {\mathcal{F}}(u^{j,-})(\xi)d\xi\right)^2=4t\left(\frac{1}{2}\int_{\xi\in\R} i\xi e^{-t|\xi|}   {\mathcal{F}}(u^{j,-})(\xi)d\xi\right)^2\nonumber\\
&=& t \left(\int_{\xi\in\R} \mathcal{F}[\partial_x G](\xi)\mathcal{F}(u^{j,-})(\xi)d\xi\right)^2\nonumber\\
&=& t \left(\int_{x\in\R}  {\partial_{x} G(x,t)}  (u^{j,-})(x)dx\right)^2\nonumber\\
&=&t\frac{1}{\pi^2} \left(\int_{x\in\R}  \frac{2xt}{(x^2+t^2)^2} (u^{j,-})(x)dx\right)^2\nonumber
\end{eqnarray}
From \rec{q},  \rec{aplus} and \rec{bminus} it follows that
\begin{equation}\label{Po1}
 \left|\int_{x\in\R}  \frac{x^2-t^2}{(x^2+t^2)^2}  (u^+)(x)dx\right|^2=  \left|\int_{x\in\R}  \frac{2xt}{(x^2+t^2)^2}  (u^-)(x)dx\right|^2
\end{equation}
In particular we get
\begin{equation}\label{Po2}
  \left|\int_{x\in\R}  \frac{x^2-t^2}{(x^2+t^2)^2}  u(x)dx\right| = \left|\int_{x\in\R}  \frac{2xt}{(x^2+t^2)^2}  u(x)dx\right|,
\end{equation}
 which achieves the proof of theorem \ref{Poho1}.\hfill$\square$
\subsection*{Preliminary computations}
We observe that
$$
\frac{x^2-t^2}{(x^2+t^2)^2}=\frac{d}{dx}\left(\frac{-x}{x^2+t^2}\right),~\frac{2xt}{(x^2+t^2)^2} =\frac{d}{dx}\left(\frac{-t}{x^2+t^2}\right).$$
 
{\bf Computation of ${\mathcal{F}}\left[\frac{-1}{x^2+1}\right]$ and ${\mathcal{F}}\left[\frac{-x}{x^2+1}\right].$}\par
\begin{eqnarray*}
{\mathcal{F}}\left[\frac{-1}{x^2+1}\right](\xi)&=&\int_{\R} \frac{-1}{x^2+1} e^{-ix\xi} d\xi\\
&=& -\pi e^{-|\xi|}.
\end{eqnarray*}
\begin{eqnarray*}
{\mathcal{F}}\left[\frac{-x}{x^2+1}\right](\xi)&=&i\frac{d}{d\xi}\left({\mathcal{F}}\left[\frac{-1}{x^2+1}\right]\right)=-i\pi\frac{\xi}{|\xi|}e^{-|\xi|}.\end{eqnarray*}

{\bf Computation of $(-\Delta)^{-1/4}[\frac{x^2-1}{(1+x^2)^2}]$ and $(-\Delta)^{-1/4}[\frac{2x}{(1+x^2)^2}]$}
\begin{eqnarray}
\label{FM+}
(-\Delta)^{-1/4}\left[\frac{x^2-1}{(1+x^2)^2}\right]&=&{\mathcal{F}}^{-1}\left[|\xi|^{-1/2}{\mathcal{F}}\left[\frac{d}{dx}\left(\frac{-x}{x^2+1}\right)\right]\right]\nonumber\\
&=&{\mathcal{F}}^{-1}[|\xi|^{-1/2}i\xi {\mathcal{F}}\left[\frac{-x}{x^2+1}\right]\nonumber\\
&=&\pi  {\mathcal{F}}^{-1}\left[|\xi|^{-1/2}\xi \frac{\xi}{|\xi|}e^{-|\xi|}\right]=\pi  {\mathcal{F}}^{-1}\left[|\xi|^{1/2}e^{-|\xi|}\right]\\
&=&\pi \frac{1}{2\pi}\int_{\R}|\xi|^{1/2}e^{-|\xi|}e^{-ix\xi} d\xi=\sqrt{\pi} \Re\left(\frac{1}{(1+ix)^{3/2}}\right)\nonumber\\
&=&\sqrt{\pi} \left(\frac{\cos(\arctan(-x))}{(1+x^2)^{3/4}}\right).\nonumber
\end{eqnarray}

\begin{eqnarray}
\label{FM-}
(-\Delta)^{-1/4}\left[\frac{2x}{(1+x^2)^2}\right]&=&{\mathcal{F}}^{-1}\left[|\xi|^{-1/2}{\mathcal{F}}\left[\frac{d}{dx}\left(\frac{-1}{x^2+1}\right)\right]\right]\nonumber\\
&=&{\mathcal{F}}^{-1}[|\xi|^{-1/2}i\xi {\mathcal{F}}\left[\frac{-1}{x^2+1}\right]\nonumber\\
&=&\pi  i{\mathcal{F}}^{-1}\left[|\xi|^{-1/2}\xi e^{-|\xi|}\right]=\pi i\frac{1}{2\pi}\int_{\R}|\xi|^{-1/2}\xi e^{-|\xi|}e^{-ix\xi} d\xi\\
&=&\sqrt{\pi}i \Im\left(\frac{1}{(1+ix)^{3/2}}\right)=\sqrt{\pi} \left(\frac{\sin(\arctan(-x))}{(1+x^2)^{3/4}}\right).\nonumber
\end{eqnarray}

Next we define the following operators:
\begin{eqnarray}
M^-[w](t)&:=&\int_{\R} \sqrt{\pi} \left(\frac{\sin(\arctan(-x))}{(1+x^2)^{3/4}}\right)w(tx) dx,\label{Phi}\\
M^+[w](t)&:=&\int_{\R} \sqrt{\pi} \left(\frac{\cos(\arctan(-x))}{(1+x^2)^{3/4}}\right)w(tx) dx,\label{Psi}.
\end{eqnarray}
\begin{Proposition}\label{M+M-}
The operators $M^{+}$ (resp. $M^-$) is an isomorphim from $L^p_+$ to $L^p_{+}$ (resp. $L^p_-$ to $L^p_{-}$)  for every $p\in (1,+\infty).$
\end{Proposition}
{\bf Proof.}
We prove the proposition for $M^+$, it is exactly the same for $M^-$.\\

{\bf Claim 1.} $M^+\colon L^p_+(\R)\to L^p_+(\R)$, for every $p\in (1,+\infty).$\\

{\bf Proof of the claim 1.}
Let $p\in (1,\infty)$ and $p'=\frac{p}{p-1}$ be the conjugate of $p$.

Let us introduce $g(x):=\frac{(1+x^2)^{\beta}}{|x|^{\alpha}}$ where $\beta=\alpha=\frac{1}{4p'},$ and let $w\in L^p(\R)$.\par 

\begin{eqnarray*}
\|M^+[w](t)\|^p_{L^p(\R^+)}&=& (\pi)^{p/2}\int_\R\left(\int_{\R} \left(\frac{\cos(\arctan(-x))}{(1+x^2)^{3/4}}\right)w(tx) dx\right)^p dt\\
&\le & (\pi)^{p/2}\left(\int_{\R}|w(y)|^p dy\right)\left(\int_{\R}\left(\frac{g(x)}{(1+x^2)^{3/4}}\right)^{p'}dx\right)^{p/p'}\left(\int_{\R}\frac{1}{|x|g^{p}(x)}\,dx\right)\\
&=& (\pi)^{p/2}\left(\int_{\R}|w(y)|^p dy\right)\\
&&\left(\int_{\R} |x|^{-1/4}(1+x^2)^{-\frac{3}{4} p'+1/4}dx\right)^{p/p'}\left(\int_{\R} |x|^{-1+\frac{p}{4p'}}(1+x^2)^{-\frac{p}{4p'}}\,dx\right)\\
&\le& C_p\|w\|^p_{L^p(\R)}.~~\hfill\Box
\end{eqnarray*}

{\bf Claim 2:  The adjoint of $M^+$  in  $L_+^2$  is $\mathcal{F}^{-1}\circ M^+ \circ \mathcal{F}^{-1} $}\\

{\bf Proof of claim 2.} 
Thanks to (\ref{FM+}), we have
$$ \pi |\xi|^{-1/2}\xi \frac{\xi}{\vert \xi\vert}e^{-|\xi|} =\mathcal{F}\left(\sqrt{\pi} \left(\frac{\cos(\arctan(-x))}{(1+x^2)^{3/4}}\right) \right).$$
Then let $w\in L^2(\R)$, for $t\not=0$,
$$M^+[w](t)=\pi \int_{\R}|\xi|^{1/2} e^{-|\xi|} {\cal{F}}[w(tx)](\xi) d\xi=\pi \int_{\R}|\xi|^{1/2} \vert t\vert^{-1}e^{-|\xi|}\hat{w}\left(\frac{\xi}{t}\right) d\xi .$$
Here we used the fact that $\hat{w}$ is real. Then, let $w,v \in L_+^2$, we have
{\begin{eqnarray*}
\langle v,M^+[w]\rangle&=&\pi \int_\R v(t)\left(\int_{\R}|\xi|^{1/2} e^{-|\xi|}\vert t\vert^{-1}\hat{w}\left(\frac{\xi}{t}\right) d\xi \right) dt\\
&=&\pi \int_\R \int_{\R} |\xi|^{1/2} e^{-|\xi|} v(t)\vert t\vert^{-1}\hat{w}\left(\frac{\xi}{t}\right) dt d\xi \\
&=&\pi \int_\R \int_{\R} |y|^{1/2} e^{-|ty|} v(t)\vert t\vert^{1/2}\hat{w}(y) dt dy\\
&=&\pi \int_\R \int_{\R}\vert x\vert^\frac{1}{2}  e^{-|x|}  |y|^{-1} v\left(\frac{x}{y}\right)\hat{w}(y) dx dy\\ 
&=& \langle M^+ \left[\mathcal{F}^{-1}(v)\right], \mathcal{F}(w)\rangle\\
&=& \left\langle \mathcal{F}^{-1}\left( M^+ \left[\mathcal{F}^{-1}(v)\right]\right), w\right\rangle
\end{eqnarray*}}

{\bf Claim 3: $Ker(M^+)=Ker(M^+)^*=\{0\}$}\\

{\bf Proof of claim 3.}
We observe that $M^+[w]=0$ if and only if the Laplace transform of $|\xi|^{-1/2}\xi \hat{w}({\xi})$ is zero. Here we use once more that $w$ is even. Since the Laplace transfom is injective, it implies that $|\xi|^{-1/2}\xi \hat{w}({\xi})=0$ a.e, namely $w= 0$ a.e. The same hold for the adjoint.\\

By combining Claim 1, 2, 3 we deduce that $M^+$ is bijective bounded linear operators from $L^2_+$ to $L^2_+$. The Open Mapping Theorem
implies that it is an isomorphism as well.
\par
By density arguments we get that $M^+$ is also an isomorphism from $L^p_+$ to $L^p_+$ for every $p\in (1,\infty)$ and therefore from
$L^{2,1}_+$ to $L^{2,1}_+$ too.
We conclude the proof of Proposition \ref{M+M-}.~\hfill $\Box$.
 
  \subsection{Pohozaev Identities for $(-\Delta)^{1/2}$ in $S^1$}
 In this section we will derive a Pohozaev formula for $C^1$, stationary $1/2$-harmonic on $S^1$, namely satisfying
 \begin{equation}
 \frac{d}{da}\int_{S^1}| (-\Delta)^{1/4}(u\circ\phi_a)|^2|_{a=0}=0
 \end{equation}
 where $\phi_a\colon S^1\to S^1$ is a family of diffeomorphisms such that 
 $\phi_0=Id.$ \par Stationary $1/2$-harmonic maps satisfy
 \begin{eqnarray}
 0&=&(-\Delta)^{1/2}(u\circ\phi_a)\cdot \frac{d}{d\theta}(u\circ\phi_a)|_{a=0}\\
 &=&(-\Delta)^{1/2}(u\circ\phi_a)\cdot \frac{du}{d\phi}\circ\phi_a\frac{d\phi_a}{d\theta}|_{a=0}~~\mbox{in}~{\cal{D}}^{\prime}\nonumber
\end{eqnarray}
If we choose for every $a\in [0,\frac{1}{2})$, $\phi_a=\frac{e^{i\theta}-a}{1-ae^{i\theta}}$ and we set 
$u_a(\theta)=u(\frac{e^{i\theta}-a}{1-ae^{i\theta}})$ (namely we consider the composition of $u$ with a
family of M\"obius transformations of the disk).\par
In this case we get 
\begin{equation}
\frac{du}{d\phi}\circ\phi_a|_{a=0}=-ie^{-i\theta}\partial_{\theta}u_0(\theta),~~\frac{d\phi_a}{d\theta}|_{a=0}=-1+e^{2i\theta}.\end{equation}
Moreover if $\phi_a$ is conform we have
\begin{equation}
(-\Delta)^{1/2}(u\circ\phi_a)=(-\Delta)^{1/2}u_a= e^{\lambda_a} ((-\Delta)^{1/2}u)\circ \phi_a,
\end{equation}
where $\lambda_a=\log(|\frac{\partial\phi_a}{\partial\theta}(\theta)|).$
Therefore
 \begin{eqnarray}\label{harms1}
  0&=&(-\Delta)^{1/2}u\circ \phi_a\cdot \frac{du}{d\phi}\circ\phi_a\frac{d\phi_a}{d\theta}|_{t=0}~~\mbox{in}~{\cal{D}}^{\prime}\nonumber\\
  &=&(-\Delta)^{1/2}u (e^{i\theta})\cdot (-ie^{-i\theta}\partial_{\theta}u_0(\theta)(-1+e^{2i\theta}))\\
  &=&(-\Delta)^{1/2}u_0 \cdot (2\sin(\theta)\partial_{\theta}u_0(\theta) ).\nonumber
\end{eqnarray}

We also  observe that
  the energy $\int_{S^1} |(-\Delta)^{1/4}u|^2\, d\theta $ is invariant with respect to the trace of M\"obius transformations of the disk and therefore for every function $w\in \dot H^{1/2}(S^1)$ we have
$$\frac{d}{da}\int_{S^1} |(-\Delta)^{1/4}w_a|^2\, d\theta=\frac{d}{da}\int_{S^1} |(-\Delta)^{1/4}w|^2\, d\theta=0.$$
In particular we get
$$2\int_{S^1}\frac{dw_a}{da}(\theta)|_{a=0}(-\Delta)^{1/2}w\, d\theta=4\int_{S^1}\sin(\theta)\partial_{\theta}w (-\Delta)^{1/2}w\, d\theta=0.$$
 
In the sequel we identify $S^1$ with $[-\pi,\pi]$. We consider the following problem
\begin{equation}\label{solfund3}
\left\{\begin{array}{cc}
\partial_t F+(-\Delta)^{1/2} F=0 & ~~\theta\in[-\pi,\pi],\; t>0\\
F(0,\theta)=\delta_0(x) & \theta\in [-\pi,\pi].
\end{array}\right.
\end{equation}
{\bf Claim:} The solution of \rec{solfund3} is given by
$$
F(\theta,t)=\frac{1}{2\pi}\sum_{n=-\infty}^{+\infty} e^{-t|n|}e^{in\theta}=\frac{e^{2t}-1}{e^{2t}-2e^{t}\cos(\theta)+1}.$$
{\bf Proof of the Claim.} 
For every $t>0$ we write
\begin{equation}\label{Fn}
F(\theta,t)=\sum_{n=-\infty}^{+\infty} F_n(t)e^{in\theta}.\end{equation}
By plugging the formula into the equation we get for every $n\in \Z$:
$$
\frac{d}{dt}F_n(t)+|n|F_n(t)=0.$$
Therefore 
$F_n(t)=C_ne^{-t|n|}$ and $F(\theta,t)=\sum_{n=-\infty}^{+\infty} C_ne^{-t|n|}e^{in\theta}.$ Since
$F(\theta,0)=\delta_0=\frac{1}{2\pi}\sum_{n=-\infty}^{+\infty} e^{in\theta}$ , namely
$  F(\theta,0)\ast f(\theta)=f(\theta)$ for every distribution $f$,  we get that $C_n=\int_{S^1} F(\theta,0)e^{-in\theta} d\theta=\frac{1}{2\pi}$ for every $n\in\Z$.
Hence 
$F$ is as in \rec{Fn} and the claim is proved.~$\Box$
\medskip
 
\begin{Theorem}{[Case on $S^1$]}\label{thms1}
Let $u\in \dot H^{1/2}_{loc}(S^1,\R^m)$ such that
\begin{equation}\label{harmequation2}
\frac{\partial u}{\partial \theta}\cdot (-\Delta)^{1/2} u=0~~\mbox{a.e in $S^1$.}\end{equation}
Then the following identity holds
\begin{equation}\label{pohos1}
 \left|\int_{S^1}u(z)\partial_t F(z) \, d\theta\right|^2= \left|\int_{S^1}u(z)\partial_{\theta}F(z) \, d\theta\right|^2
\end{equation}
\end{Theorem}
{\bf Proof of Theorem \ref{thms1}}
 
 Now we write
$u(\theta)=\sum_{n\in\Z} u_ne^{in\theta}$. Then the following equalities hold:
\begin{eqnarray}\label{harms2}
\sin(\theta)\partial_{\theta}u &=&\frac{1}{2}\sum_{n\in\Z}\left(e^{i\theta}-e^{-i\theta}\right) nu_ne^{in\theta}\nonumber\\
&=& \frac{1}{2}\sum_{n\in\Z}e^{in\theta}((n-1)u_{n-1}-(n+1)u_{n+1})\nonumber\\ 
\overline{(-\Delta)^{1/2}u}&=&\sum_{n\in\Z} |n|\bar{u}_ne^{-in\theta}\nonumber\\
\sin(\theta)\partial_{\theta}u\overline{(-\Delta)^{1/2}u}&=&\sum_{n,m\in\Z}e^{i(n-m)\theta}|m|[(n-1)u_{n-1}-(n+1)u_{n+1}]\bar{u}_m.
\end{eqnarray}
By combining  \rec{harms1} and \rec{harms2} we get 
\begin{eqnarray}\label{harms3}
0&=&\sum_{n,m\in\Z}|m|  e^{-t|n-m|}\left[(n-1)u_{n-1}u_{n-1}-(n+1)u_{n=1}\right]\bar {u}_{m}\nonumber\\
&=&\sum_{n,m\in\Z}|m|nu_n\bar{u}_{m}\left[e^{t|n-m+1|}-e^{-t|n-m-1|}\right].\end{eqnarray}
We    first  symmetrize \rec{harms3}. We get first
\begin{equation} \label{harms4}
0=- \sum_{n,m\in\Z}|n|mu_m\bar{u}_{n}\left[e^{t|n-m+1|}-e^{-t|n-m-1|}\right].
\end{equation}
We sum \rec{harms3} and \rec{harms4}
  
  \begin{equation} \label{harms5}
0=\sum_{n,m\in\Z}|  \left[e^{t|n-m+1|}-e^{-t|n-m-1|}\right] [|m|n u_n\bar{u}_{m}-|n|m\bar {u}_nu_m].
\end{equation}
 
For every $n\in\Z$ we set 
\begin{eqnarray*}
u_n&=&a_n+ib_n\\
a_n^+&=&\frac{a_n+a_{-n}}{2},~~~a_n^-=\frac{a_n-a_{-n}}{2}\\
b_n^+&=&\frac{b_n+b_{-n}}{2},~~b_n^-=\frac{b_n-b_{-n}}{2}.
\end{eqnarray*}
We observe that since the components of $u$ are  real we have $a^-_n=0$ and $b^+_n=0$, $a_n^+=u^+_n$, $b_n^-=i u_n^-.$
We take the real part of \rec{harms5} and get

 \begin{eqnarray} 
0&=& 
\sum_{n,m\in\Z}(|m|n-|n|m) \left[e^{-t|n-m+1|}-e^{-t|n-m-1|}\right] [a_ma_n+b_nb_m]\nonumber\\
&=& \sum_{n,m\in\Z}(|m|n-|n|m) \left[e^{-t|n-m+1|}-e^{-t|n-m-1|}\right] [a^+_ma^+_n] \label{harms6}\\
&&+ \sum_{n,m\in\Z}(|m|n-|n|m) \left[e^{-t|n-m+1|}-e^{-t|n-m-1|}\right] [b^-_mb^-_n].\label{harms7}
\end{eqnarray}
 {\bf Estimate of \rec{harms6}}
  \begin{eqnarray} \label{harms8}
&&  \sum_{n,m\in\Z}(|m|n-|n|m) \left[e^{-t|n-m+1|}-e^{-t|n-m-1|}\right] [a^+_ma^+_n] \nonumber\\&=&
 \sum_{n>0,m<0}(-2mn) \left[e^{-t|n-m+1|}-e^{-t|n-m-1|}\right] [a^+_ma^+_n] \nonumber\\
 &&~~+\sum_{n<0,m>0}(2mn) \left[e^{-t|n-m+1|}-e^{-t|n-m-1|}\right] [a^+_ma^+_n] \nonumber\\
 &=& \sum_{n>0,m<0}(-2mn) \left[e^{-t|n-m+1|}-e^{-t|n-m-1|}\right] [a^+_ma^+_n] \nonumber\\
&&~~+\sum_{n>0,m<0}(2(-m)(-n)) \left[e^{-t|-n+m+1|}-e^{-t|-n+m-1|}\right] [a^+_{-m}a^+_{-n}n]\nonumber\\
&=& -4\sum_{n>0,m<0}(mn) \left[e^{-t|n-m+1|}-e^{-t|n-m-1|}\right] [a^+_ma^+_n] \nonumber\\
&=& 4\sum_{n>0,m>0}(-m)n \left[e^{-t|n+m+1|}-e^{-t|n+m-1|}\right] [a^+_{-m}a^+_n] 
\nonumber\\
&=& -4\sum_{n>0,m>0}mn \left[e^{-t(n+m+1)}-e^{-t(n+m-1)|}\right] [a^+_{m}a^+_n]
\nonumber\\
&=& -4\left(e^{-t}-e^{t}\right)\left(\sum_{n>0}ne^{-tn}  a^+_n\right)^2=8\sinh(t)\left(\sum_{n\in\Z}|n|e^{-t|n|}  u^+_n\right)^2\nonumber\\
&=&8\sinh(t)\left(\sum_{n\in\Z}(-\Delta)^{1/2}F_n  u^+_n\right)^2\nonumber\\
&=&\frac{1}{2\pi}8\sinh(t)\left(\int_0^{2\pi}(-\Delta)^{1/2}F  u^+ d\theta\right)^2=\frac{1}{2\pi}8\sinh(t)\left(\int_0^{2\pi}{\partial_t F}  u^+ d\theta\right)^2.
\end{eqnarray}

{\bf Estimate of \rec{harms7}}
  \begin{eqnarray} \label{harms9}
&&  \sum_{n,m\in\Z}(|m|n-|n|m) \left[e^{-t|n-m+1|}-e^{-t|n-m-1|}\right] [b^-_mb^-_n] \nonumber\\&=&
 \sum_{n>0,m<0}(-2mn) \left[e^{-t|n-m+1|}-e^{-t|n-m-1|}\right] [b^-_mb^-_n] \nonumber\\
 &&~~+\sum_{n<0,m>0}(2mn) \left[e^{-t|n-m+1|}-e^{-t|n-m-1|}\right] [b^-_mb^-_n] \nonumber\\
 &=& \sum_{n>0,m<0}(-2mn) \left[e^{-t|n-m+1|}-e^{-t|n-m-1|}\right] [b^-_mb^-_n] \nonumber\\
&&~~+\sum_{n>0,m<0}(2(-m)(-n)) \left[e^{-t|-n+m+1|}-e^{-t|-n+m-1|}\right] [b^-_{-m}b^-_{-n}n]\nonumber\\
&=& -4\sum_{n>0,m<0}(mn) \left[e^{-t|n-m+1|}-e^{-t|n-m-1|}\right] [b^-_mb^-_n] \nonumber\\
&=& 4\sum_{n>0,m>0}(-m)n \left[e^{-t|n+m+1|}-e^{-t|n+m-1|}\right] [b^-_{-m}b^-_n] 
\nonumber\\
&=& 4\sum_{n>0,m>0}mn \left[e^{-t(n+m+1)}-e^{-t(n+m-1)|}\right] [b^-_{m}b^-_n]
\nonumber\\
&=& 4\left(e^{-t}-e^{t}\right)\left|\sum_{n>0}ne^{-tn}  b^-_n\right|^2=-8\sinh(t)\left|\sum_{n\in\Z}ine^{-t|n|}  u^-_n\right|^2\nonumber\\
&=&-8\sinh(t)\left|\sum_{n\in\Z}(\partial_{\theta} F_n ) u^-_n\right|^2\nonumber\\
&=&-\frac{1}{2\pi}8\sinh(t)\left|\int_0^{2\pi}(\partial_{\theta} F)  u^- d\theta\right|^2.
\end{eqnarray}

From \rec{harms6} and \rec{harms8} and \rec{harms9} it follows
\begin{equation}\label{pohozaev2}
\left|{ \int_0^{2\pi}{\partial_t F}} (t,\theta) u(\theta) d\theta\right| =\left| \int_0^{2\pi}{{\partial_{\theta} F}} (t,\theta) u(\theta) d\theta\right|
 \end{equation}
where
$${\partial_t F} (t,\theta)=-2e^t\frac{e^{2t}\cos(\theta)-2e^{t}+\cos(\theta)}{(e^{2t}-2e^{t}\cos(\theta)+1)^2}$$
and
$$
{\partial_{\theta} F} (t,\theta) =-2e^{t}\frac{ \sin(\theta)(e^{2t}-1)}{(e^{2t}-2e^{t}\cos(\theta)+1)^2}.$$
    \subsection{Pohozaev Identities for the Laplacian in $\R^2$}
    In this section we derive a Pohozaev identity in $2D$ which analogous to that found in $1D$. Instead of integrating on balls,  the strategy 
    is to multiply by the fundamental solution of the heat equation. A similar idea has been preformed in \cite{Stru} to study the heat flow.\par
The solution of
\begin{equation}\label{solfundlap}
\left\{\begin{array}{cc}
\partial_t G+(-\Delta) G=0 & t>0\\
G(0,x)=\delta_{x_0} & t=0\,.
\end{array}\right.
\end{equation}
is given by $G(x,t)=(4\pi t)^{-1/2}e^{-\frac{|x-x_0|^2}{4t}}.$
\begin{Theorem}{[Case on $\R^2$]}\label{Poho0}
Let $u\in \C^2(\R^2,\R^m)$ such that
\begin{equation}\label{harmequation0}
\left\langle \frac{\partial u}{\partial x_i},\Delta u\right\rangle=0~~\mbox{a.e in $\R^2$}
\end{equation}
$i=1,2.$
Assume that
\begin{equation}\label{cond1}
\int_{\R}|u-u_0|dx <+\infty, ~~\int_{\R}|\nabla u(x)|dx<+\infty.
\end{equation}
Then for all $x_0\in\R^2$ and $t>0$ the  following identity holds
\begin{equation}\label{idpohozaevR2}
  \iint_{R^2}e^{-{\frac{|x-x_0|^2}{4t}}}|x-x_0|^2\left|\frac{\partial u}{\partial \nu} \right|^2dx=  \iint_{R^2}e^{-{\frac{|x-x_0|^2}{4t}}}\left|\frac{\partial u}{\partial \theta}\right|^2 dx.
  \end{equation}
\end{Theorem}
{\bf Proof.}
We multiply the equation \rec{harmequation0} by $x_i e^{-{\frac{|x-x_0|^2}{4t}}}$ and we integrate:
\begin{eqnarray}\label{estpoh0}
0&=&\sum_{k,i=1}^2\iint_{\R^2}  e^{-{\frac{|x-x_0|^2}{4t}}}(x_i -{x_0}_i)\frac {\partial u}{\partial x_i}\frac{\partial^2 u}{\partial x_k^2}dx\nonumber\\
&=&-\frac{1}{2t}\iint_{\R^2}  e^{-{\frac{|x-x_0|^2}{4t}}}|x -x_0|^2\left|\frac{\partial u}{\partial\nu}\right|^2 dx]-\iint_{\R^2}e^{-{\frac{|x-x_0|^2}{4t}}}|\nabla u|^2 dx\nonumber\\
&-&\frac{1}{2}\sum_{i=1}^2\iint_{\R^2}  e^{-{\frac{|x-x_0|^2}{4t}}}(x -x_0)_i\frac{\partial}{\partial x_i} |\nabla u|^2 dx\nonumber\\
&=&-\frac{1}{2t}\iint_{\R^2}  e^{-{\frac{|x-x_0|^2}{4t}}}|x -x_0|^2\left|\frac{\partial u}{\partial\nu}\right|^2 dx +\frac{1}{4t}\iint_{\R^2}  e^{-{\frac{|x-x_0|^2}{4t}}}|x -x_0|^2\left|\nabla u\right|^2 dx.
\end{eqnarray}
Since $\nabla u=(\frac{\partial u}{\partial\nu},|x-x_0|^{-1}\frac{\partial u}{\partial\theta})$, from \rec{estpoh0} we get the identity
\begin{equation} 
  \iint_{R^2}e^{-{\frac{|x-x_0|^2}{4t}}}|x-x_0|^2\left|\frac{\partial u}{\partial \nu} \right|^2dx=  \iint_{R^2}e^{-{\frac{|x-x_0|^2}{4t}}}\left|\frac{\partial u}{\partial \theta}\right|^2 dx
  \end{equation}
  and we conclude.~\hfill$\Box$\par
  \medskip
 \begin{Remark}
We obtain an analogous identity to \rec{idpohozaevR2} if  we multiply the equation \rec{harmequation0} by $X_i e^{-{\frac{|x-x_0|^2}{4t}}}$ where
  $X_1+i X_2\colon \C\to \C$ is a holomorpic function.
  By using the Cauchy-Riemman differential equations we get
  \begin{equation}
  2  \iint_{R^2}e^{-{\frac{|x-x_0|^2}{4t}}}|x-x_0| \left\langle \frac{\partial u}{\partial \nu}, \frac{\partial u}{\partial X}\right\rangle dx=   \iint_{R^2}e^{-{\frac{|x-x_0|^2}{4t}}}\left\langle x-x_0, X\right\rangle |\nabla u|^2 dx.
  \end{equation}
  \end{Remark}

  \section{ Compactness and Quantization of horizontal $1/2$ harmonic maps}\label{cq}
 
 The proof  of the first part of Theorem \ref{th-I.2} is exactly the same of that of Lemma 2.3 in \cite{DL2} and we omit it. \par
 As far as the quantization issue  is concerned the proof goes as that of Theorem 1.1 in \cite{DL2} (namely the decomposition of the $\R$  into converging regions, bubbles domains and neck-regions) once we perform the analysis
 of the neck-region. As we have already mentioned in the Introduction global uniform $L^{2,1}$ estimates in degenerating annuli are not anymore available as in the case of $1/2$-harmonic maps with values into a sphere.
 Therefore we have to perform a subtle analysis of ${(-\Delta)}^{1/4}u$ where $u\in {\frak{H}}^{1/2} ({\R})$ is  a  horizontal $1/2$-harmonic map in an annular domain. To this purpose we will make use of the
 Pohozaev-type formulae we have discovered in Section \ref{Poh}.

  We first show that if $u\in {\frak{H}}^{1/2} ({\R})$ is  a  horizontal $1/2$-harmonic map with  ${(-\Delta)}^{1/2}u\in L^{1}(\R)$ then ${\cal{R}}[ ( P^N (-\Delta)^{1/4}u)\in L^{2,1}(\R).$
 
  Next we show that if $u\in \dot H^{1/2}(\R,{\cal{N}})$ is a weak harmonic map, then ${(-\Delta)}^{1/2}u\in L^{1}(\R)$. \par
\begin{Proposition}\label{prL1}
Let  $u\in \dot H^{1/2}(\R,{\cal{N}})$ be a weak harmonic maps. Then ${(-\Delta)}^{1/2}u\in L^{1}(\R)$.  
\end{Proposition}
{\bf Proof of Proposition \ref{prL1}.}
{\bf Step 1.}  We prove that ${(-\Delta)}^{1/2}u\in L^{1}(\R)$.\par
This follows from the fact that for all $\xi,\eta\in {\cal{N}}$ we have 
$$P^N(\xi)\cdot (\xi-\eta)=O(|\xi-\eta|^2).$$
If  $u\in \dot H^{1/2}(\R,{\cal{N}})$ is a weak harmonic map, then
$P^N(u){(-\Delta)}^{1/2}u ={(-\Delta)}^{1/2}u$ in  ${\cal{D}}^\prime(\R)$. Therefore
\begin{eqnarray}\label{estL1}
\int_{\R}|{(-\Delta)}^{1/2}u(x)| dx&=& \int_{\R}|P_N(u(x)){(-\Delta)}^{1/2}u(x)| dx\\
&\le &  C\int_{\R}\int_{\R}\frac{|u(x)-u(y)|^2}{|x-y|^2} dx dy \nonumber\\
&=& C\|u\|^2_{\dot{H}^{1/2}(\R)}<+\infty.\nonumber~\Box
\end{eqnarray}
\begin{Proposition}\label{prL2}
If $u\in {\frak{H}}^{1/2} ({\R})$ is  a  horizontal $1/2$-harmonic map with  ${(-\Delta)}^{1/2}u\in L^{1}(\R)$ then ${\cal{R}}[ ( P^N (-\Delta)^{1/4}u)]\in L^{2,1}(\R).$
\end{Proposition}
{\bf Proof of Proposition \ref{prL2}.}
 Since $u\in {\frak{H}}^{1/2} ({\R})$, we have $P_N(u)\nabla u=P_N(u){\cal{R}}(-\Delta)^{1/2} u=0.$
 The result follows from the fact that $P_N(u)(-\Delta)^{1/4} u$ satisfies the following structure equation
 \begin{equation}\label{eqstruct}
 (-\Delta)^{1/4} ( P^N (-\Delta)^{1/4}u)= S(P^N,u)-{\cal{R}} [((-\Delta)^{1/4} P^N) ({\cal{R}}(-\Delta)^{1/4} u)]\,.
\end{equation}
 We deduce  in particular that ${\cal{R}} [((-\Delta)^{1/4} P^N) ({\cal{R}}(-\Delta)^{1/4} u)]\in L^{1}(\R)$.\par Since $((-\Delta)^{1/4} P^N) ({\cal{R}}(-\Delta)^{1/4} u)\in L^1(\R)$ it follows that
 $((-\Delta)^{1/4} P^N) ({\cal{R}}(-\Delta)^{1/4} u)\in {\mathcal{H}}^1$.\footnote{See section \ref{secdef} for the definition.}The Hardy Space can be also characterized as the space of function in $f\in L^1$ such that
${\cal{R}}[f]\in L^1({\R}) $.\par
Hence \rec{eqstruct}  implies that 
$ (-\Delta)^{1/4}{\cal{R}}[ ( P^N (-\Delta)^{1/4}u)]\in {\mathcal{H}}^1(\R)$ and hence ${\cal{R}}[ ( P^N (-\Delta)^{1/4}u)] \in L^{2,1}(\R).$ We conclude the proof.~\hfill$\Box$\par
\bigskip

 {\bf Proof of Theorem \ref{th-I.3}}.
 We have already remarked that the sequence  $$v_k=(P^k_T(-\Delta)^{1/4}u_k,{\cal{R}}[P^k_N(-\Delta)^{1/4} u_k])$$ satisfies  a system of the form
  \begin{equation}\label{modelsystemk}
(-\Delta)^{1/4} v_k=\Omega_k v_k+\Omega^1_k  v+{\cal{Z}}(Q_k,v_k) \end{equation}
where $v_k\in L^2(\R)$, $Q_k\in  \dot H^{1/2}(\R)$,  $\Omega_k\in L^2(\R,so(m))$, $\Omega_k^1\in L^{2,1}(\R)$, $g_k\in L^1(\R)$ and ${\cal{Z}}\colon  \dot H^{1/2}(\R)\times L^2(\R)\to {\cal{H}}^1(\R)$ with
\begin{equation}
 \|\Omega_k\|_{L^2(\R)}+\|\Omega^1_k\|_{L^{2,1}(\R)}+\|Q_k\|_{\dot H^{1/2}(\R)} )\le C\|v_k\|_{L^2}
\end{equation}
If $\delta>0$ is small enough we can apply Proposition \ref{decomposition}.

\begin{equation}\label{estdeltauring}
(P^k_T+ {\mathcal R} {P^k_N} ) (-\Delta)^{1/4} u_k(x)=\11_{A_{r_k,R_k}}[A^{-1}(x)\overrightarrow{c_{r_k}}\frac{1}{|x|^{1/2}}]+h_k(x) +g_k(x)\end{equation}
with
$A^{-1}_k\in L^{\infty}\cap \dot H^{1/2}(\R,M_m({\R}))),$ $h_k\in L^{2,1}(A_{\Lambda r_k, (2\Lambda)^{-1}R_k}),$ (for every $\Lambda>2$ such that $\Lambda r_k<(2\Lambda)^{-1}R_k$),  $g_k\in L^2(\R)$ with ${\rm supp} g_k\subset (A_{r,R})^c.$ Moreover  we have 
\be\label{estck}
\displaystyle \overrightarrow{c_{r_k}}=O\left(\left(\log(\frac{R_k}{2\Lambda^2 r_k})\right)^{-1/2}\right), ~~\mbox{as}~k\to +\infty, \Lambda\to +\infty. \ee
Therefore
\begin{eqnarray} 
P^k_T (-\Delta)^{1/4} u_k(x)&=&P^k_T\left(\frac{   \11_{A_{r_k,R_k}{A_k^{-1}(x)}\overrightarrow{c_{r_k}}}}{|x|^{1/2}}\right)+P^k_T h(x)+P^k_T g(x)\label{PTu}\\
P^k_N (-\Delta)^{1/4} u_k(x)&=&P^k_N\left({\mathcal{R}}\left[\frac{ \11_{A_{r_k,R_k}}{A^{-1}(x)}\overrightarrow{c_{r_k}}}{|x|^{1/2}}\right]\right)-P^k_N{\mathcal{R}}[ h(x)]-P^k_N{\mathcal{R}}[  g(x)].\label{PNu}
\end{eqnarray}
Observe that $P^k_N{\mathcal{R}} g_k(x)\in L^{2,1}(A_{\Lambda r_k, (2\Lambda)^{-1}R_k})$ since ${\rm supp}( g_k)\subset (A_{r_k,R_k})^c.$
By combining \rec{PTu} and \rec{PNu} we get
\begin{equation}
  (-\Delta)^{1/4}u_k(x)=P^k_T\left(\frac{  \11_{A_{r_k,R_k}}{A_k^{-1}(x)}\overrightarrow{c_{r_k}}}{|x|^{1/2}}\right) +\bar h_k \ee
   with $\bar h_k=P^k_T h_k(x)+P^k_T(x)+P^k_N (-\Delta)^{1/4} u_k(x)  $  which is in $ L^{2,1}(A_{\Lambda r_k, (2\Lambda)^{-1}R_k}).$
 
Next
we set
$a_k(x):= P^k_T{A_k^{-1}}$  and we denote for simplicity by $a_k^+$ and $a_k^-$ respectively the symmetric and antisymmetric parts of $a_k$.
Since $A_k^{-1}$ and $P^k_T$ are in $\dot H^{1/2}(\R)$, they verify the following estimate (see \cite{Tartar}):
   \begin{eqnarray}
   \left\|\frac{A_k^{-1}(x)-A_k^{-1}(-x)}{|x|^{1/2}}\right\|_{L^2}&\le& C\|A_k^{-1}\|_{\dot H^{1/2}}\\
     \left\|\frac{P^k_T(x)-P^k_T(-x)}{|x|^{1/2}}\right\|_{L^2}&\le& C\|P^k_T\|_{\dot H^{1/2}}\label{tarpt}\\
  \end{eqnarray}
By using the fact that  $A_k^{-1}$ and $P^k_N,P^k_T$ are also in $L^{\infty}$ we get that
$$\left\| \frac{a_k^{-}\overrightarrow{c_{r_k}}}{|x|^{1/2}}\right\|_{L^2(\R)}\le C|\overrightarrow{c_{r_k}}|.$$
 Therefore we can write
 \begin{equation}
   (-\Delta)^{1/4}u_k(x)=\11_{A_{r_k,R_k}} \frac{a_k^{+}\overrightarrow{c_{r_k}}}{|x|^{1/2}}  +\bar h_k+\tilde g_k
    \end{equation}
where $\bar h_k\in L^{2,1}(A_{\Lambda r_k, (2\Lambda)^{-1}R_k}),$ with
$$\limsup_{\Lambda\to\infty}\limsup_{k\to\infty } \|  \bar h_k\|_{L^{2,1}(B(0,(2\Lambda)^{-1}R_k)\setminus B(0,\Lambda r_k)} <+\infty,$$ and 
$\tilde g_k=\11_{A_{r_k,R_k}} \frac{a_k^{-}\overrightarrow{c_{r_k}}}{|x|^{1/2}}\in L^{2}(\R)$, with $\|\tilde g_k\|_{L^2}\le |\overrightarrow{c_{r_k}}|.$
We can conclude.~~\hfill$\Box$
\par
\bigskip
{\bf Proof of Theorem \ref{th-I.5}}\par
 Under the current hypotheses Theorem \ref{th-I.3} and Lemma \ref{neck2} it follows :\par   
\begin{equation}\label{linftyestbis}
\limsup_{\Lambda\to\infty}\limsup_{k\to\infty } \| (-\Delta)^{1/4} u_k\|_{L^{2,\infty}(B(0,(2\Lambda)^{-1}R_k)\setminus B(0,\Lambda r_k)} =0\end{equation}
   and
     \begin{equation}\label{devel3}
   (-\Delta)^{1/4}u_k(x)=\11_{A_{r_k,R_k}} \frac{a_k^{+}\overrightarrow{c_{r_k}}}{|x|^{1/2}}  +\bar h_k+\tilde g_k,
    \end{equation}
where $\bar h_k\in L^{2,1}(A_{\Lambda r_k, (2\Lambda)^{-1}R_k}),$ with
$$\limsup_{\Lambda\to\infty}\limsup_{k\to\infty } \|  \bar h_k\|_{L^{2,1}(B(0,(2\Lambda)^{-1}R_k)\setminus B(0,\Lambda r_k)} <+\infty,$$
 $\tilde g_k\in L^{2}(\R)$ and $\|\tilde g_k\|_{L^2}\le C |\overrightarrow{c_{r}}|.$
   By combining \rec{linftyestbis} and \rec{devel3} we get
   \begin{equation}\label{estasymi}
\limsup_{\Lambda\to\infty}\limsup_{k\to\infty } \| ((-\Delta)^{1/4} u_k)^-\|_{L^{2}(B(0,(2\Lambda)^{-1}R_k)\setminus B(0,\Lambda r_k)} =0\end{equation}
In order to establish  a link between the symmetric and antisymmetric part of $(-\Delta)^{1/4} u_k$ we make use of the formula \ref{idpohozaevR}
than can be rewritten as
\begin{equation}\label{refPoh}
 \left(\int_{\R} \sqrt{\pi} \left(\frac{\sin(\arctan(-x))}{(1+x^2)^{3/4}}\right)w(tx) dx\right)^2=\left(\int_{\R} \sqrt{\pi} \left(\frac{\cos(\arctan(-x))}{(1+x^2)^{3/4}}\right)w(tx) dx\right)^2
\end{equation}
Now we plug into \rec{refPoh} the function $w_k(x)=\11_{A_{\Lambda r_k, (2\Lambda)^{-1}R_k}} (-\Delta)^{1/4}u_k(x)$. We observe that
 $$\int_{\R} \sqrt{\pi} \left(\frac{\sin(\arctan(-x))}{(1+x^2)^{3/4}}\right)\11_{A_{r_k,R_k}}(tx) \frac{a_k^{+}\overrightarrow{c_{r_k}}}{|tx|^{1/2}} dx=0.$$
Therefore we have that 
\begin{eqnarray*}
&&\int_{\R} \sqrt{\pi} \left(\frac{\cos(\arctan(-x))}{(1+x^2)^{3/4}}\right)\11_{A_{r_k,R_k}}(tx) \frac{a_k^{+}\overrightarrow{c_{r_k}}}{|tx|^{1/2}} dx\\
&=&
\int_{\R} \sqrt{\pi} \left(\frac{\sin(\arctan(-x))}{(1+x^2)^{3/4}}\right)\tilde{w}_k(tx) dx - \int_{\R} \sqrt{\pi} \left(\frac{\cos(\arctan(-x))}{(1+x^2)^{3/4}}\right)\tilde{w}_k(tx) dx.
\end{eqnarray*}
where
$$\tilde w_k(x)=w_k(x)-\11_{A_{r_k,R_k}}(x) \frac{a_k^{+}\overrightarrow{c_{r_k}}}{|x|^{1/2}}=\bar h_k+\tilde g_k.$$
Next we use the fact that the operators 
\begin{eqnarray}
M^-[w](t)&:=&\int_{\R} \sqrt{\pi} \left(\frac{\sin(\arctan(-x))}{(1+x^2)^{3/4}}\right)w(tx) dx,\label{M+}\\
M^+[w](t)&:=&\int_{\R} \sqrt{\pi} \left(\frac{\cos(\arctan(-x))}{(1+x^2)^{3/4}}\right)w(tx) dx,\label{M-}.
\end{eqnarray}
are isomorphism from $L^{p}_+$ to $ L^{p}_+$ (resp. $L^{p}_-$ to $ L^{p}_-$) for every $p>1$ and from  $L^{2,1}_+$ to $ L^{2,1}_+$ (resp.  $L^{2,1}_-$ to $ L^{2,1}_-$ ) and we deduce that
$$
\11_{A_{r_k,R_k}}(x) \frac{a_k^{+}\overrightarrow{c_{r_k}}}{|x|^{1/2}}=\varphi_k+\psi_k.$$
with  $\varphi_k\in L^{2,1}(A_{\Lambda r_k, (2\Lambda)^{-1}R_k}),~\psi_k\in L^{2}(\R)$ and $\|\psi_k\|_{L^2}\le C|\overrightarrow{c_{r}}|.$
\par
Hence
\begin{eqnarray}
&&\limsup_{\Lambda\to\infty}\limsup_{k\to+\infty}\|(-\Delta)^{1/4}u_k-\psi_k\|_{L^2(B(0,(2\Lambda)^{-1}R_k)\setminus B(0,\Lambda_k r_r))}
\\&&\le \limsup_{\Lambda\to\infty}\limsup_{k\to+\infty}\|(-\Delta)^{1/4}u_k-\psi_k\|_{L^{2,\infty}(B(0,(2\Lambda)^{-1}R_k)\setminus B(0,\Lambda_k r_r))}
\nonumber \\
&&\cdot  \limsup_{\Lambda\to\infty}\limsup_{k\to+\infty}\|(-\Delta)^{1/4}u_k-\psi_k\|_{L^{2,1}(B(0,(2\Lambda)^{-1}R_k)\setminus B(0,\Lambda_k r_r))}=0.
\end{eqnarray}
Therefore since $\limsup_{k\to+\infty}\|\psi_k\|_{L^2}=0$ we deduce that
$$
\limsup_{\Lambda\to\infty}\limsup_{k\to+\infty}\|(-\Delta)^{1/4}u_k\|_{L^2(B(0,(2\Lambda)^{-1}R_k)\setminus B(0,\Lambda_k r_r))}.$$

We can conclude the proof of Theorem \ref{th-I.5}.~~\hfill $\Box$
\section{Counter-example}
The aim of this section is to construct a sequence of solutions  of a Schr\"odinger type equation with antisymmetric potential whose energy  is not quantized. In particular, we are going to build  a sequence of solutions whose energy of the potential goes to zero in the neck but not the energy of the solutions.\\

We consider
$u(t)= 1 $ on $[-1,1]$ and $\frac{1}{\vert t\vert^{\frac{1}{2}}}$ elsewhere and  $v(t)= \frac{1}{(1+t^2)^{\frac{3}{8}}}$.

\begin{Lemma} As $t\rightarrow +\infty$, we have
 $$(-\Delta)^\frac{1}{4} u(t)  =O\left( \frac{1}{t^\frac{3}{2}}\right)$$
 and
  $$(-\Delta)^\frac{1}{4} v(t)  =O\left( \frac{1}{t^\frac{5}{4}}\right).$$
\end{Lemma}

{\it Proof :}\\
\begin{equation}
\begin{split}
(-\Delta)^\frac{1}{4} u(t) &=\int_{\R} \frac{u(t)-u(s)}{\vert t-s\vert^\frac{3}{2}} \, ds \\
&=\int_{-\infty}^{-1} \frac{u(t)-u(s)}{\vert t-s\vert^\frac{3}{2}} \, ds + \int_{-1}^{1} \frac{u(t)-u(s)}{\vert t-s\vert^\frac{3}{2}} \, ds+ \int_{1}^{+\infty} \frac{u(t)-u(s)}{\vert t-s\vert^\frac{3}{2}} \, ds\\
& =I_1 + I_2+I_3.
\end{split}
\end{equation} 
Let $ t >1$

\begin{equation}
\begin{split}
I_2 &=  \int_{-1}^1 \frac{t^{-\frac{1}{2}}-1}{\vert t-s\vert^\frac{3}{2}} \, ds \\
&=\frac{t^{-\frac{1}{2}} -1}{ t^\frac{3}{2}}\int_{-1}^{1} \frac{1}{\vert1 -\frac{s}{t}\vert^\frac{3}{2}} \, ds\\
&=\frac{t^{-\frac{1}{2}} -1}{ t^\frac{1}{2}}\int_{-\frac{1}{t}}^{\frac{1}{t}} \frac{1}{\vert1 - u\vert^\frac{3}{2}} \, du\\
&=\frac{t^{-\frac{1}{2}} -1}{ t^\frac{1}{2}} \left[ 2(1-u)^{-\frac{1}{2}} \right]^{\frac{1}{t}}_{-\frac{1}{t}}
=O\left( \frac{1}{t^\frac{3}{2}} \right)
\end{split}
\end{equation}

\begin{equation}
\begin{split}
I_1 &= \int_{-\infty}^{-1} \frac{t^{-\frac{1}{2}}-(-s)^{-\frac{1}{2}}}{\vert t-s\vert^\frac{3}{2}} \, ds \\
&=\frac{1}{t^2}\int_{-\infty}^{-1} \frac{1-(-\frac{s}{t})^{-\frac{1}{2}}}{\vert1 -\frac{s}{t}\vert^\frac{3}{2}} \, ds\\
&=\frac{1}{t} \int^{+\infty}_{\frac{1}{t}} \frac{1-u^{-\frac{1}{2}}}{(1+u)^\frac{3}{2}} \, du 
\end{split}
\end{equation} 

\begin{equation}
\begin{split}I_3 &= \int_{1}^{+\infty} \frac{t^{-\frac{1}{2}}-s^{-\frac{1}{2}}}{\vert t-s\vert^\frac{3}{2}} \, ds \\ 
&=\frac{1}{t^2}\int_{1}^{+\infty} \frac{1-(\frac{s}{t})^{-\frac{1}{2}}}{\vert1 -\frac{s}{t}\vert^\frac{3}{2}} \, ds\\
&=\frac{1}{t} \int^{+\infty}_{\frac{1}{t}} \frac{1-u^{-\frac{1}{2}}}{\vert 1-u\vert^\frac{3}{2}} \, du 
 \end{split}
\end{equation} 
We easily check that $\displaystyle \frac{1-u^{-\frac{1}{2}}}{\vert 1-u\vert^\frac{3}{2}} \underset{u\rightarrow 1}{\sim} (1-u)^{-\frac{1}{2}}$ and then the last integral is well defined.\\

But changing the variable $u$ into $\frac{1}{v}$ into $I_1$ and $I_2$ we observe that 

$$  \int^{+\infty}_{0} \frac{1-u^{-\frac{1}{2}}}{(1+u)^\frac{3}{2}} \, du =  \int_{0}^{+\infty} \frac{1-u^{-\frac{1}{2}}}{\vert 1-u\vert^\frac{3}{2}} \, du=0$$

Which implies that 
$$I_1=\frac{-1}{t} \int_0^{\frac{1}{t}} \frac{1-u^{-\frac{1}{2}}}{(1+u)^\frac{3}{2}} \, du =  \underset{t\rightarrow +\infty}{O} \left(\frac{1}{t^\frac{3}{2}}\right)$$
and

$$I_3=\frac{-1}{t} \int_0^{\frac{1}{t}} \frac{1-u^{-\frac{1}{2}}}{\vert 1-u\vert^\frac{3}{2}} \, du   =  \underset{t\rightarrow +\infty}{O} \left(\frac{1}{t^\frac{3}{2}}\right).$$
Which proves the lemma for $u$ when $t>1$. Of course we have the same result for $t<-1$ by symmetry.\\

 \begin{equation}
\begin{split}
(-\Delta)^\frac{1}{4} v(t) &=\int_{\R} \frac{u(t)-u(s)}{\vert t-s\vert^\frac{3}{2}} \, ds=  \frac{1}{(1+t^2)^{\frac{3}{8}}} \int_{\R} \frac{(1+s^2)^{\frac{3}{8}} -(1+t^2)^{\frac{3}{8} }}{ (1+s^2)^{\frac{3}{8}}\vert t-s\vert^\frac{3}{2}} \, ds\\
&= \frac{1}{(1+t^2)^{\frac{3}{8}}} \int_{-\infty}^{-1} \frac{(1+s^2)^{\frac{3}{8}} -(1+t^2)^{\frac{3}{8} }}{ (1+s^2)^{\frac{3}{8}}\vert t-s\vert^\frac{3}{2}} \, ds + \frac{1}{(1+t^2)^{\frac{3}{8}}} \int_{-1}^1 \frac{(1+s^2)^{\frac{3}{8}} -(1+t^2)^{\frac{3}{8} }}{ (1+s^2)^{\frac{3}{8}}\vert t-s\vert^\frac{3}{2}} \, ds \\
&+ \frac{1}{(1+t^2)^{\frac{3}{8}}} \int_{1}^{+\infty} \frac{(1+s^2)^{\frac{3}{8}} -(1+t^2)^{\frac{3}{8} }}{ (1+s^2)^{\frac{3}{8}}\vert t-s\vert^\frac{3}{2}} \, ds
\\& =I_1 + I_2+I_3.
\end{split}
\end{equation} 

We remark  that
$$ \frac{(1+s^2)^{\frac{3}{8}} -(1+t^2)^{\frac{3}{8} }}{\vert t-s\vert^\frac{3}{4}} \hbox{ is uniformly bounded.}$$
Then, let $t>1$,
\begin{equation}
\begin{split}
I_2 &= O\left( \frac{1}{\vert t \vert^{\frac{3}{4}}} \int_{-1}^1 \frac{1}{ (1+s^2)^{\frac{3}{8}}\vert t-s\vert^\frac{3}{4}} \, ds\right)\\
&= O\left( \frac{1}{\vert t \vert^{\frac{3}{4}}} \int_{-1}^1 \frac{1}{ \vert t-s\vert^\frac{3}{4}} \, ds\right)\\
&= O\left( \frac{1}{\vert t \vert^{\frac{5}{4}}} \right) 
\end{split}
\end{equation} 
And

\begin{equation}
\begin{split}
I_1 &= O\left( \frac{1}{\vert t \vert^{\frac{3}{4}}} \int_{-\infty}^{-1} \frac{1}{ (1+s^2)^{\frac{3}{8}}\vert t-s\vert^\frac{3}{4}} \, ds\right)\\
&= O\left( \frac{1}{\vert t \vert^{\frac{3}{4}}} \int_{-\infty}^{-1} \frac{1}{ (-s)^\frac{3}{4} \vert t-s\vert^\frac{3}{4}} \, ds\right)\\
&= O\left( \frac{1}{\vert t \vert^{\frac{3}{4}}} \int_{-\infty}^{-1} \frac{1}{ (-s)^\frac{3}{4} \vert t-s\vert^\frac{3}{4}} \, ds\right)\\
&= O\left( \frac{1}{\vert t \vert^\frac{5}{4} } \int_{-\infty}^{-1/t} \frac{1}{ (-u)^\frac{3}{4} \vert 1-u\vert^\frac{3}{4}} \, ds\right)\\
&= O\left( \frac{1}{\vert t \vert^{\frac{5}{4}}} \right) 
\end{split}
\end{equation} 
The same estimate work for $I3$ and $t<-1$, which achieved the proof.\hfill$\square$

Then we set $\omega =\frac{(-\Delta)^\frac{1}{4} u}v$ and $\omega_1= \frac{(-\Delta)^\frac{1}{4} v +\omega u}{u}$, which gives

$$(-\Delta)^\frac{1}{4} U=\Omega U +\Omega_1 U$$
with
$$U=\left(\begin{array}{c}u \\v\end{array}\right),$$

$$\Omega = \left(\begin{array}{cc}0 & \omega \\-\omega & 0\end{array}\right)$$
and
$$\Omega_1=\left(\begin{array}{cc}0 & 0 \\ \omega_1 & 0\end{array}\right).$$

Finally, we set  $U_n(t)=c_n U(nt)$, $\Omega_n= \sqrt{n} \Omega(nt)$, and ${\Omega_1}_n=\sqrt{n}  \Omega_1(nt)$. Where
$$c_n =\frac{1}{\Vert u(nt)\Vert_2}= \left( \frac{n}{\ln\left(\frac{n+\sqrt{1+n^2}}{-n+\sqrt{1+n^2}}\right)}\right)^\frac{1}{2} \sim \left( \frac{n}{\ln(n)}\right)^\frac{1}{2} $$

We easily check that

$$\Vert U_n \Vert_2 \sim1$$

$$\lim_{R\rightarrow +\infty} \lim_{n\rightarrow +\infty} \Vert \Omega_n \Vert_{L^2(B(0,1/R) \setminus B(0,R/n))} +  \Vert {\Omega_1}_n \Vert_{L^{2,1} (B(0,1/R) \setminus B(0,R/n))} =0$$
and 
$$(-\Delta)^\frac{1}{4} U_n = \Omega_n U_n +{\Omega_1}_n U_n$$ 

which prove that there is no quantification of the energy, since $\Omega_n$  and  ${\Omega_1}_n$ have no energy in the neck region despite $U_n$ get some. Of course such an example doesn't satisfies a Pohozaev identity, since it is symmetric. It is a "local" example since $U_n$ bounded in $L^2$ only on $[-1,1]$, it would be interesting to construct a example define on the whole $\R$. It may be possible by starting on the circle and then projecting on to $\R$. Indeed the main idea here is to take  a modification of the Green function of the $(-\Delta)^{\frac{1}{4}}$. Hence playing this game on $S^1$ and projecting then $\R$ should provide an example with good decreasing at infinity.

\end{document}